\documentclass[10pt]{amsart}
\usepackage{amssymb,amsmath,amsthm}
\usepackage{mathrsfs,dsfont,a4wide}
\theoremstyle{plain}
\newtheorem{theorem}{Theorem}[section]
\newtheorem{lemma}[theorem]{Lemma}
\newtheorem{proposition}[theorem]{Proposition}
\newtheorem{corollary}[theorem]{Corollary}

\theoremstyle{definition}
\theoremstyle{remark}
\numberwithin{equation}{section}

\newcommand{\as}{{\mathcal A}}
\newcommand{\hs}{{\mathcal H}}

\newcommand{\gs}{{\mathcal G}}
\newcommand{\fs}{{\mathcal F}}
\newcommand{\leb}{{\mathcal L}}

\newcommand{\ms}{{\mathcal M}}

\newcommand{\eub}{{\mathcal E}}

\newcommand{\Es}{{\mathcal E}}

\newcommand{\dst}{{\mathcal D}}

\newcommand{\R}{{\mathbb R}}
\newcommand{\N}{{\mathbb N}}

\newcommand{\epi}[1]{{\rm epi}(#1)}
\newcommand{\dom}[1]{{\rm dom}(#1)}
\newcommand{\1}{{\mathds 1}}

\newcommand{\Om}{\Omega}
\newcommand{\Omb}{\overline{\Omega}}
\newcommand{\Omt}{\tilde{\Omega}}
\newcommand{\Omp}{\Omega'}

\newcommand{\weakst}{\stackrel{\ast}{\rightharpoonup}}

\newcommand{\weak}{\rightharpoonup}
\newcommand{\wlystar}{$\text{weakly}^*$\;}
\newcommand{\wstar}{$\text{weak}^*$\;}

\newcommand{\eps}{\varepsilon}
\newcommand{\tsub}{\,\tilde{\subseteq}\,}
\newcommand{\Sg}[2]{S^{#1}(#2)}

\input psfig.sty

\newcommand{\res}{\mathop{\hbox{\vrule height 7pt width .5pt depth 0pt
\vrule height .5pt width 6pt depth 0pt}}\nolimits}

\title
[Size effects on quasistatic growth of fractures]
{Size effects on quasistatic growth of fractures}
\author[A. Giacomini]
{Alessandro Giacomini}
\address[Alessandro Giacomini]{S.I.S.S.A., Via Beirut 2-4, 34014, Trieste,
Italy}
\email[A. Giacomini]{giacomin@sissa.it}
\begin{document}
\vskip .2truecm
\begin{abstract}
\small{
We perform an analysis of the size effect for quasistatic growth
of fractures in linearly isotropic elastic bodies under antiplanar shear.
In the framework of the variational model proposed by G.A. Francfort
and J.-J. Marigo in \cite{FM},
we prove that if the size of the body tends to infinity,
and even if the surface energy is of cohesive form,
under suitable boundary displacements the fracture propagates
following the
Griffith's functional.
}
\end{abstract}
\maketitle
{\small \tableofcontents}

\section{Introduction}
\label{intr}
A well known fact in fracture mechanics is that
{\it ductility} is also influenced by the size of the structure,
and in particular
the structure tends to become brittle if its size increases
(see for example \cite{Ca}, and
references therein).
The aim of this paper is to capture this fact for the problem of 
quasistatic growth of fractures
in linearly elastic bodies in the framework of the variational 
theory of crack propagation
formulated by Francfort and Marigo in \cite{FM}.
\par
The model proposed in \cite{FM} is inspired to classical Griffith's 
criterion
and determines the evolution of the fracture through a competition
between volume and
surface energies. Let us illustrate it and the variant we investigate
in the case of {\it generalized antiplanar shear}.
\par
Let $\Om \subseteq \R^N$ be open, bounded, and with Lipschitz boundary.
A fracture $\Gamma \subseteq \Omb$
is any rectifiable set, and a displacement $u$ is any function
defined almost everywhere in $\Om$ whose
set of discontinuities $S(u)$ is contained in $\Gamma$
(we will make precise the functional
setting later). The total energy of the configuration $(u,\Gamma)$
is given by
\begin{equation}
\label{basicenintr}
\int_{\Om \setminus \Gamma} |\nabla u|^2\,dx+\hs^{N-1}(\Gamma).
\end{equation}
The first term in \eqref{basicenintr} implies that we assume to
apply linearized elasticity
in the unbroken part of $\Om$. The second term can be considered
as the work done to create $\Gamma$.
\par
As suggested in \cite{FM}, more general fracture energies can be
considered in \eqref{basicenintr},
especially those of Barenblatt's type \cite{Ba},
and here we consider energies of the form
\begin{equation}
\label{barenbintr}
\int_\Gamma \varphi(|[u]|(x)) \,d\hs^{N-1}(x),
\end{equation}
where $[u](x):=u^+(x)-u^-(x)$ is the difference 
of the traces of $u$
on both sides of $\Gamma$, and
$\varphi:[0,+\infty[ \to [0,+\infty[$ 
(which depends on the material)
is such that $\varphi(0)=0$.
In order to get a physical interpretation of
\eqref{barenbintr},
let us
set $\sigma:=\varphi'$: we interpret 
$\sigma(|[u]|(x))$ as density
of force
in $x$ that act between the two lips of the crack $\Gamma$ whose
displacement are
$u^+(x)$ and $u^-(x)$ respectively. Typically $\sigma$
is decreasing, and
$\sigma(s)=0$ for
$s \ge \bar s$: this means that the interaction between
the two lips of
the fracture decreases
as the opening increases, and disappear when the opening
is greater than a
critical length $\bar s$.
As a consequence, $\varphi$ is increasing and concave,
and $\varphi(s)$ is
constant for $s \ge \bar s$.
We will then consider $\varphi$ increasing, concave,
with $\varphi(0)=0$,
$c=\varphi'(0)<+\infty$, and
$\lim_{s \to +\infty} \varphi(s)=1$. We can interpret
$$
\int_\Gamma \varphi(|[u]|(x))\,d\hs^{N-1}(x)
$$
as the work made to create $\Gamma$ with an opening given
by $[u]$. Assuming linearized elasticity to hold in
$\Om \setminus \Gamma$, we consider a total energy of the form
\begin{equation}
\label{deftotenintr}
\|\nabla u\|^2+\int_\Gamma \varphi(|[u]|) \,d\hs^{N-1},
\end{equation}
where $\|\cdot\|$ denotes the $L^2$ norm.
The problem of irreversible quasistatic growth of fractures in the
cohesive case
can be addressed through a {\it time discretization process}
in analogy to what proposed in \cite{FM}
for the energy \eqref{basicenintr}.
\par
Let $g(t)$ be a time dependent boundary displacement defined on
$\partial_D \Om \subseteq \partial \Om$
with $t \in [0,T]$. Let $\delta>0$ and let
$I_\delta:=\{0=t^\delta_0<t^\delta_1<\dots<t^\delta_{N_\delta}=T\}$
be a subdivision of $[0,T]$ with
$\max (t^\delta_{i+1}-t^\delta_i)<\delta$, and let
$g^\delta_i:=g(t^\delta_i)$. Let us consider a preexisting
crack configuration $(\bar \Gamma,\bar \psi)$,
with $\bar \psi$ a positive function
on $\bar \Gamma$: this means that
$\int_{\bar \Gamma}\varphi(\bar \psi)\,d\hs^{N-1}$
is the work done on the fracture $\bar \Gamma$ before
we apply the boundary displacement $g(t)$.
At time $t=0$ we consider $u^\delta_0$ as a minimum of
\begin{equation}
\label{def1intr}
\|\nabla u\|^2+ \int_{\Sg{g(0)}{u} \cup \bar \Gamma}
\varphi(|[u]| \lor \bar \psi)\,d\hs^{N-1}.
\end{equation}
Here $\Sg{g(0)}{u}:=
S(u) \cup \{x \in \partial_D \Om\,:\,u(x) \not= g(0)(x)\}$,
and
for all $x \in \partial_D \Om$ we consider
$[u](x):=g(x)-\tilde{u}(x)$, where $\tilde{u}$ is the
trace of $u$ on $\partial \Om$. Moreover we intend that
$\bar \psi=0$ outside $\bar \Gamma$.
We define the fracture $\Gamma^\delta_0$ at time $t=0$ as
$\Sg{g(0)}{u^\delta_0} \cup \bar \Gamma$. We also set
$\psi^\delta_0:=\bar \psi \lor |[u^\delta_0]|$ on
$\Gamma^\delta_0$.
The presence of $\Sg{g(0)}{u}$ in \eqref{def1intr}
indicates that the points at which the boundary displacement
is not attained are considered as a part of the fracture.
Notice moreover that problem \eqref{def1intr} takes into
account an {\it irreversibility condition} assumption
in the growth of the fracture.
Indeed, while on $\Sg{g(0)}{u} \setminus \bar \Gamma$
the surface energy which comes in minimization
of \eqref{def1intr} is exactly as in \eqref{barenbintr}, on
$\Sg{g(0)}{u} \cap \bar \Gamma$ the surface energy
involved takes into account the previous work made on $\bar \Gamma$.
The surface energy is of the form
of \eqref{barenbintr} only if $|[u]|>\bar \psi$, that is only if the
opening is increased. If $|[u]| \le \bar \psi$
no energy is gained, that is displacements of this form along the fracture
are in a sense surface energy free. Notice
finally that the irreversibility condition involves only the modulus of $[u]$:
this is an assumption
which is reasonable since we are considering only antiplanar displacement.
Clearly more complex irreversibility conditions can be formulated, involving
for example a partial release of energy: the one we study is the first
straightforward extension of the irreversibility condition given in \cite{FM}
for the energy \eqref{basicenintr}.
\par
Supposing to have constructed $\Gamma^\delta_i$ and $\psi^\delta_i$ at time
$t^\delta_i$, we consider
a minimum $u^\delta_{i+1}$ of the problem
\begin{equation}
\label{def2intr}
\|\nabla u\|^2+ \int_{\Sg{g(t^\delta_{i+1})}{u} \cup  \Gamma^\delta_i}
\varphi(|[u]| \lor
\bar \psi^\delta_i)\,d\hs^{N-1}
\end{equation}
and define $\Gamma^\delta_{i+1}:=
\Gamma^\delta_i \cup \Sg{g(t^\delta_{i+1})}{u^\delta_{i+1}}$ and
$\psi^\delta_{i+1}:=\psi^\delta_i \lor |[u^\delta_{i+1}]|$ on
$\Gamma^\delta_{i+1}$. As noted previously,
problem \eqref{def2intr} involves an irreversibility condition:
the surface energy density on $\Gamma^\delta_i$ increases
only where $|[u]|> \psi^\delta_i$.
\par
The {\it discrete in time evolution}
of the fracture relative to the boundary datum $g(t)$,
the subdivision $I_\delta$,
the initial crack configuration $(\bar \Gamma,\bar \psi)$ is given by
$\{(u^\delta_i, \Gamma^\delta_i,\psi^\delta_i)\,:\,i=0,\dots,N_\delta\}$.
\par
The {\it irreversible quasistatic evolution} of fracture relative
to the boundary datum $g(t)$ and the initial crack configuration
$(\bar \Gamma,\bar \psi)$
is obtained as a limit for $\delta \to 0$ of
$(u^\delta(t),\Gamma^\delta(t),\psi^\delta(t))$,
where
$u^\delta(t):=u^\delta_i$, $\Gamma^\delta(t):=\Gamma^\delta_i$ and
$\psi^\delta(t):=\psi^\delta_i$ for $t^\delta_i \le t<t^\delta_{i+1}$.
\par
This program has been studied in detail in several papers in the case
$\varphi \equiv 1$, that is for energy of the
form \eqref{basicenintr}. A first mathematical formulation
has been given in
\cite{DMT},
where the authors consider the case of dimension $N=2$ and
fractures which
are compact and with a uniform bound
on the number of connected components. This analysis has been
extended to the
case of plane elasticity in \cite{Ch}.
In \cite{FL} the authors consider the general dimension $N$,
and remove the bound on the number of
the connected components of the fractures: the key point is to
introduce
a weak formulation of
the problem considering displacements in the space $SBV$
(see Section \ref{prel}).
Finally in \cite{DMFT} the authors treat the case of finite elasticity
not restricted
to antiplanar shear, with volume energy depending on the full gradient
under suitable growth condition,
and in presence of volume forces and traction forces:
the appropriate functional space for the displacements is now $GSBV$
(see for example \cite{AFP}
for a precise definition).
\par
In all these papers (\cite{DMT},\cite{Ch},\cite{FL},\cite{DMFT}),
the analysis of the limit reveals
three basic properties (irreversibility, minimality and
nondissipativity,
see Theorem \ref{qse}) which are taken as definition of irreversible
quasistatic growth of brittle fractures:
the time discretization procedure is considered as a privileged way to
get an existence result.
\par
In the case of energy \eqref{deftotenintr}, several difficulties arise
in the analysis of the
discrete in time evolution, and in the analysis as $\delta \to 0$.
In Section \ref{discrevolsec}, we prove
that the functional space we need for the step by step minimization
is the space of functions with
bounded variation $BV$ (see Section \ref{prel}): moreover we prove
that a relaxed version
of \eqref{deftotenintr} has to be employed, namely
\begin{equation}
\label{relenintr}
\int_{\Om} f(\nabla u)\,dx+\int_{\Gamma}
\varphi(|[u]| \lor \psi)\,d\hs^{N-1}+a|D^cu|(\Om),
\end{equation}
where $a=\varphi'(0)$, $f$ is defined in \eqref{bulkrelax2},
and $D^cu$ indicates the Cantorian part
of the derivative of $u$. An existence result for
{\it discrete in time evolution} in this context
of $BV$ space is given in Proposition \ref{discrevol}.
\par
The analysis for $\delta \to 0$ presents several difficulties,
the main one being the stability of the minimality
property of the discrete in time evolutions.
The main purpose of this paper is to prove that these difficulties
disappear as the size of the reference configuration
increases, thank to the fact that the body response tends to become
more and more brittle in spite of the presence of
cohesive forces on the fractures given by the function $\varphi$.
More precisely we consider a crack configuration
$(\bar \Gamma,\bar \psi)$ in $\Om$ with
$\hs^{N-1}(\bar \Gamma)<+\infty$, and prove this
fact for the discrete evolutions in $\Om_h:=h\Om$
with preexisting crack configuration $(\bar \Gamma_h,\bar \psi_h)$
of the form $\bar \Gamma_h:=h\bar \Gamma$ and $\bar \psi_h(x):=
\bar \psi (\frac{x}{h})$, under suitable boundary displacements.
The idea is to rescale displacements and fractures to the 
fixed configuration $\Om$, and take advantage from the form of the 
problem in this new setting.
The boundary displacements on $\partial_D \Om_h:=
h\partial_D \Om$ will be taken of the form
$$
g_h(t,x):=h^\alpha g\left( t, \frac{x}{h} \right),
\quad\quad
g \in AC([0,T];H^1(\Om)),
\quad\quad
\|g(t)\|_\infty \le C,
\quad\quad
t \in [0,T],\; x \in \Om_h,
$$
where $\alpha>0$ and $C>0$. We indicate by
$(u^{\delta,h}(t), \Gamma^{\delta,h}(t),\psi^{\delta,h}(t))$
the piecewise constant interpolation
of the discrete in time evolution of fracture in $\Om_h$ relative
to the boundary
displacement $g_h$ and the
preexisting crack configuration $(\bar \Gamma_h,\bar \psi_h)$.
Let us moreover set for every $t \in [0,T]$
$$
\Es^{\delta,h}(t):=\int_{h\Om}f(\nabla u^{\delta,h}(t))\,dx+
\int_{\Gamma^{\delta,h}(t)}\varphi(\psi^{\delta,h}(t))\,d\hs^{N-1}
+a|D^c u^{\delta,h}|(\Om_h).
$$
In the case $\alpha=\frac{1}{2}$, we make the following rescaling
$$
v^{\delta,h}(t,x):=\frac{1}{\sqrt{h}}u^{\delta,h}(t,hx),
\quad\quad
K^{\delta,h}(t):=\frac{1}{h}\Gamma^{\delta,h}(t),
\quad\quad x \in \Om.
$$
The main result of the paper is the following (see Theorem \ref{main1}
for a more precise statement).

\begin{theorem}
\label{main1intr}
If $\delta \to 0$ and $h \to +\infty$,
there exists a quasistatic evolution of brittle fractures
$\{t \to (v(t), K(t))\}$
in $\Om$ relative to the preexisting crack $\bar \Gamma$ and boundary
displacement $g$
in the sense of \cite{FL} (see Theorem \ref{qse}) such that for all
$t \in [0,T]$ we have
\begin{equation}
\label{gradconv1intr}
\nabla v^{\delta,h}(t) \weak \nabla v(t) \quad\quad 
\text{weakly in }L^1(\Om;\R^N),
\end{equation}
Moreover for all $t \in [0,T]$ we have
\begin{equation}
\label{totconv1intr}
\frac{1}{h^{N-1}} \eub^{\delta,h}(t) \to \|\nabla v(t)\|^2
+\hs^{N-1}(K(t));
\end{equation}
in particular $h^{-N+1}|D^cu^{\delta,h}(t)|(\Om_h) \to 0$,
\begin{equation}
\label{bulkconv1intr}
\frac{1}{h^{N-1}}\int_{\Om_h} f(\nabla u^{\delta,h}(t))\,dx \to
\|\nabla v(t)\|^2,
\end{equation}
and
\begin{equation}
\label{fracconv1intr}
\frac{1}{h^{N-1}}\int_{\Gamma^{\delta,h}(t)} \varphi(\psi^{\delta,h}(t))
\,d\hs^{N-1}
\to \hs^{N-1}(K(t)).
\end{equation}
\end{theorem}

Theorem \ref{main1intr} proves that as the size of the 
reference configuration
increases, the response of the body in the problem of 
quasistatic growth
of fractures tends to become brittle,
so that energy \eqref{basicenintr} can be considered.
Moreover we have convergence results for the volume and surface 
energies involved.
\par
The particular value $\alpha=\frac{1}{2}$ comes out
because a problem of quasistatic evolution
has been considered. In fact if we consider
an infinite plane with a crack-segment of length $l$ and subject
to a uniform stress $\sigma$ at infinity,
following Griffith's theory the crack propagates quasistatically
if $\sigma=\frac{K_{IC}}{\sqrt{\pi l}}$,
where $K_{IC}$ is the critical {\it stress intensity factor}.
So if the crack has length $hl$, the stress rescale as
$\frac{1}{\sqrt{h}}$. This is precisely what we are prescribing
in the case $\alpha=\frac{1}{2}$: in fact
the stress that intuitively we prescribe at the boundary can be
reconstructed from $\nabla u_h$ and
rescales precisely as $\frac{1}{\sqrt{h}}$.
\par
For the proof of Theorem \ref{main1intr}, the first step is to 
recognize that
$(v^{\delta,h}(t),K^{\delta,h}(t),\psi^{\delta,h}(t))$ is a
discrete in time evolution relative to the boundary displacement 
$g$ and the preexisting crack configuration
$(\bar \Gamma,\bar \psi)$ for a total energy of the form
$$
\int_{\Om} f_h(\nabla u)\,dx+\int_{\Gamma} \varphi_h(|[u]| \lor \psi)
\,d\hs^{N-1}+a\sqrt{h}|D^cu|(\Om),
$$
where $\varphi_h(s) \nearrow 1$ for all $s \in [0,+\infty[$, 
and $f_h(\xi) \nearrow |\xi|^2$ for all $\xi \in \R^N$.
From the fact that $\varphi_h \nearrow 1$ we recognize that 
the structure tends to become brittle.
Bound on total energy for the discrete in time evolution is available, 
so that compactness in the space $BV$ can be applied:
it turns out that the limits of the displacements are of
class $SBV$ with gradient in $L^2(\Om;\R^N)$. Limits for the fractures 
are constructed through a $\Gamma${-}convergence
procedure (see Lemma \ref{gammat}). Now the main point is to recover 
the minimality property (see point
(c) of Theorem \ref{qse})
$$
\|\nabla v(t)\|^2 \le \|\nabla v\|^2 +\hs^{N-1}(\Sg{g(t)}{v} 
\setminus K(t)),
\quad\quad
v \in SBV(\Om)
$$
from the minimality property of 
$(v^{\delta,h}(t),K^{\delta,h}(t),\psi^{\delta,h}(t))$.
This is done in Lemma \ref{mintlem} by means of a refined version
of the Transfer of Jump Lemma of \cite{FL}: the main difference 
here is that we have to deal with
$BV$ functions and we have
to transfer the jump on the part of $K^{\delta,h}(t)$ where
$\psi^{\delta,h}(t)$ is greater than a given small constant.
\par
We also consider the cases $\alpha \in ]0,\frac{1}{2}[$ and 
$\alpha >\frac{1}{2}$.
It turns out that in the case $\alpha \in ]0,\frac{1}{2}[$,
the body is not solicited enough to make
the preexisting crack $\bar \Gamma_h$ propagate, and $\Om_h$
tends to behave elastically in the complement of
$\bar \Gamma_h$: more precisely we prove that
(Theorem \ref{main2}) in the case $\bar \psi >\eps>0$, setting
\begin{equation}
\label{rescalpha}
v^{\delta,h}(t,x):=\frac{1}{h^\alpha}u^{\delta,h}(t,hx),
\end{equation}
for all $t \in [0,T]$ we have that $v^{\delta,h}(t)$ converges
to the displacement of the elastic problem in
$\Om \setminus \bar \Gamma$ under boundary displacement
given by $g(t)$.
\par
In the case $\alpha>\frac{1}{2}$ we have that
the preexisting fracture
$\bar \Gamma_h$ tends to propagate
brutally toward {\it rupture}: in fact in Theorem \ref{main3}
we prove that
$v^{\delta,h}(0)$ given by \eqref{rescalpha}
converges to a piecewise constant function $v$ in $\Om$, so that
$\bar \Gamma \cup
\Sg{g(0)}{v}$ disconnects $\Om$. This phenomenon is a
consequence of
the variational approach based on
the search for global minimizers: as the size of $\Om_h$
increases,
fractures carry an energy of
order $h^{N-1}$, while non rigid displacements carry an
energy of greater order: in this way fracture
is preferred to deformation.
\par
The paper is organized as follows: in Section \ref{prel} 
we recall
some basic definitions and introduce
the functional setting for the problem.
In Section \ref{discrevolsec} we deal with
the problem of
discrete
in time evolutions for fractures
in the cohesive case.
The main theorems are listed in Section \ref{mainressec},
and Sections \ref{main1sec},
\ref{main2sec} and \ref{main3sec} are devoted to their proofs.
In Section \ref{relres}
we prove a relaxation
result which is used in the problem of discrete
in time evolution of fractures,
while in Section \ref{auxsec} we prove
some auxiliary results employed in the study of the asymptotic
behavior of the evolutions.

\section{Preliminaries}
\label{prel}
In this section we state the notation and prove some preliminary 
results employed in
the rest of the paper.

\vskip20pt
{\bf Basic notation.}
We will employ the following basic notation:
\begin{itemize}
\item[-] $\Om$ is an open and bounded subset of $\R^N$ with 
Lipschitz boundary;
\item[-] $\partial_D \Om$ is a subset of $\partial \Om$ open 
in the relative topology;
\item[-] $\hs^{N-1}$ is the $(N-1)${-}dimensional
Hausdorff measure;
\item[-] we say that $A \tsub B$ if $A \subseteq B$ 
up to a set of
$\hs^{N-1}${-}measure zero;
\item[-] $\Gamma \subseteq \Om$ is rectifiable if
there exists a sequence of $C^1$ manifolds $(M_i)_{i \in \N}$
such that $\Gamma \tsub \cup_i M_i$;
\item[-] for all $A \subseteq \R^N$, $|A|$ denotes the 
Lebesgue measure of $A$;
\item[-] for all $A \subseteq \R^N$, $\1_{A}$ denotes 
the characteristic function of $A$;
\item[-] if $\mu$ is a measure on $\R^N$ and $A$ is a 
Borel subset of $\R^N$,
$\mu \res A$ denotes the restriction of $\mu$ to $A$, i.e.
$(\mu \res A)(B):=\mu(B \cap A)$ for all Borel sets 
$B \subseteq \R^N$;
\item[-] $\|u\|_{\infty}$ and $\|u\|$ denote the sup-norm and
the $L^2$ norm of $u$ respectively;
\item[-] if $u,g \in BV(\Om;\R^m)$, $\Sg{g}{u}:=
S(u) \cup \{x \in \partial_D \Om\,:\,u(x)\not=g(x)\}$;
\item[-] if $a,b \in \R$, $a \lor b:=\max\{a,b\}$ and 
$a \wedge b:=\min\{a,b\}$.
\end{itemize}

\vskip20pt
{\bf Functions of bounded variation.}
For the general theory of functions of bounded variation, 
we refer to \cite{AFP}; here
we recall some basic definitions and theorems we need in the sequel.
We say that $u \in BV(A)$ if $u \in L^1(A)$, and its distributional
derivative $Du$ is a bounded vector-valued Radon measure on $A$.
In this case it turns out that the set $S(u)$ of points $x \in A$ 
which are not
Lebesgue points of $u$ is rectifiable, that is
there exists a sequence of $C^1$ manifolds $(M_i)_{i \in \N}$ 
such that $S(u) \subseteq \cup_i M_i$ up to a set of
$\hs^{N-1}${-}measure zero.
As a consequence $S(u)$ admits a normal $\nu_u(x)$
at $\hs^{N-1}${-}a.e. $x \in S(u)$. Moreover for 
$\hs^{N-1}$ a.e. $x \in S(u)$,
there exist $u^+(x), u^-(x) \in \R$ such that
$$
\lim\limits_{r \to 0}
\frac{1}{|B^\pm_r(x)|} \int_{B^\pm_r(x)} |u(y)-u^\pm(x)|\,dy=0,
$$
where $B^\pm_r(x):=\left\{ y \in B_r(x)\,:\,
(y-x)\cdot \nu_u(x) \gtrless 0 \right\}$,
and $B_r(x)$ is the ball with center $x$ and radius $r$.
It turns out that $Du$ can be represented as
$$
Du(A)= \int_A \nabla u(x) \,dx+ 
\int_{A \cap S(u)} (u^+(x)-u^-(x)) \nu_u(x)
\,d\hs^{N-1}(x)+D^cu(A),
$$
where $\nabla u$ denotes the approximate gradient of $u$
and $D^cu$ is the Cantor part of $Du$. $BV(A)$ is a Banach space
with respect to the norm $\|u\|_{BV(A)}:=\|u\|_{L^1(A)}+|Du|(A)$.
\par
We will often use the following result: if $A$ is bounded
and Lipschitz, and if $(u_k)_{k \in \N}$ is a
bounded sequence in $BV(A)$, then there exists a subsequence
$(u_{k_h})_{h \in \N}$ and $u \in
BV(A)$ such that
\begin{align}
\label{weak*conv}
u_{k_h} \to u  \quad\quad&\text{strongly in }L^1(A), \\
\nonumber
Du_{k_h} \weakst Du  \quad\quad&
\text{\wlystar in the sense of measures}.
\end{align}
We say that $u_k \weakst u$ \wlystar in $BV(A)$ if 
\eqref{weak*conv} holds.
\par
We say that $u \in SBV(A)$ if $u \in BV(A)$ and $D^cu=0$.
The space $SBV(A)$ is called the space of
{\it special functions of bounded variation}. Note that if 
$u \in SBV(A)$, then
the singular part of $Du$ is concentrated on $S(u)$.
\par
The space $SBV$ is very useful when dealing with variational 
problems involving
volume and surface energies because of the 
following compactness and lower
semicontinuity result due to L.Ambrosio (\cite{A1}, \cite{A3}).

\begin{theorem}
\label{SBVcompact}
Let $A$ be an open and bounded subset of $\,\R^N$, and 
let $(u_k)_{k \in \N}$ be a sequence in
$SBV(A)$. Assume that there exists $q>1$ and 
$c \in [0;+\infty[$ such that
$$
\int_A |\nabla u_k|^q \,dx+ 
\hs^{N-1}(S(u_k))+ ||u_k||_\infty \le c
$$
for every $k \in \N$. Then there exists a subsequence
$(u_{k_h})_{h \in \N}$ and a function
$u \in SBV(A)$ such that
\begin{align}
\label{sbvconv}
\nonumber
u_{k_h} \to u \quad {strongly \; in}\; L^1(A), \\
\nabla u_{k_h} \weak \nabla u \quad {weakly \; in}\; 
L^1(A;\R^N), \\
\nonumber
\hs^{N-1}(S(u)) \le \liminf_h \hs^{N-1}(S(u_{k_h})).
\end{align}
\end{theorem}

In the rest of the paper, we will say that $u_k \weak u$
weakly in $SBV(A)$ if
$u_k$ and $u$ satisfy \eqref{sbvconv}.
It will also be useful the following fact which can be 
derived from
Ambrosio's Theorem: if $u_k \weak u$ weakly in $SBV(A)$ 
and if
$\hs^{N-1} \res S(u_k) \weakst \mu$ \wlystar in the
sense of measures, then $\hs^{N-1} \res S(u) \le \mu$ 
as measures.
\par
Finally in the context of fracture problems we will use
the following notation: if $A$ is Lipschitz, and if
$\partial_D A \subseteq \partial A$, then for all 
$u,g \in BV(A)$ we set
\begin{equation}
\label{sgjump}
\Sg{g}{u}:=S(u) \cup 
\{x \in \partial_D A\,:\,u(x) \not= g(x)\},
\end{equation}
where the inequality on $\partial_D A$ is intended in 
the sense of traces. Moreover, we set
for all $x \in S(u)$
$$
[u](x):=u^+(x)-u^-(x),
$$
and for all $x \in \partial_D A$ we set $[u](x):=u(x)-g(x)$, 
where the traces of $u$ and $g$ on $\partial A$
are used.

\vskip20pt
{\bf Quasi-static evolution of brittle fractures.}
Let $\Om$ be an open bounded subset of $\R^N$ with Lipschitz 
boundary, and let
$\partial_D \Om$ be a subset of $\partial \Om$ open in the 
relative topology.
Let $g:[0,T] \to H^1(\Om)$ be absolutely continuous (see \cite{Br1}
for a precise definition); we indicate
the gradient of $g$ at time $t$ by $\nabla g(t)$, and the time
derivative of $g$ at time $t$ by $\dot{g}(t)$. For $u \in SBV(\Om)$,
let $\Sg{g(t)}{u}$ be defined
as in \eqref{sgjump}, and for every $A,B \subseteq \R^N$, let
$A \tsub B$ mean $A \subseteq B$ up to a set of $\hs^{N-1}${-}measure
zero. The main result of \cite{FL} is
the following theorem.

\begin{theorem}
\label{qse}
Let $\bar \Gamma$ be a rectifiable set in $\Om \cup \partial_D \Om$ 
such that
$\hs^{N-1}(\bar \Gamma)<+\infty$. There exists 
$\{(u(t),\Gamma(t))\,:\,t \in [0,T]\}$
with $\Gamma(t) \tsub \Om \cup \partial_D \Om$ rectifiable and
$u(t) \in SBV(\Om)$ with $\Sg{g(t)}{u(t)} \tsub \Gamma(t)$
such that:
\begin{itemize}
\item[(a)] $\bar \Gamma \tsub \Gamma(s) \tsub \Gamma(t)$ for 
all $0 \le s \le t \le T$;
\item[{}]
\item[(b)] $u(0)$ minimizes
$$
\|\nabla v\|^2+\hs^{N-1}( \Sg{g(0)}{v} \setminus \bar \Gamma)
$$
among all $v \in SBV(\Om)$;
\item[{}]
\item[(c)]
for $t \in ]0,T]$, $u(t)$ minimizes
$$
\|\nabla v\|^2+\hs^{N-1} \left( \Sg{g(t)}{v} \setminus \Gamma(t)
\right)
$$
among all $v \in SBV(\Om)$.
\end{itemize}
Furthermore, the total energy
$$
\Es(t):= \|\nabla u(t)\|^2 +\hs^{N-1}( \Gamma(t))
$$
is absolutely continuous and satisfies
\begin{equation}
\label{enevolprel}
\Es(t)=\Es(0)+2 \int_0^t \int_\Om \nabla u(\tau)
\nabla \dot{g}(\tau) \,dx \,d\tau
\end{equation}
for every $t \in [0,T]$.
\end{theorem}

Condition $(a)$ stands for the {\it irreversibility} of the crack 
propagation, conditions
$(b)$ and $(c)$ are {\it minimality} conditions, while 
\eqref{enevolprel} stands
for the {\it nondissipativity} of the process.

\vskip20pt
\noindent
{\bf $\Gamma${-}convergence.}
Let us recall the definition and some basic properties
of De Giorgi's
{\it $\Gamma$-convergence} in metric spaces.
We refer the reader to \cite{dm} for an exhaustive
treatment of this subject.
Let $(X,d)$ be a metric space. We say that a sequence
$F_h:X\to [-\infty ,+\infty
]$ $\Gamma $-converges to $F:X\to [-\infty ,+\infty ]$
(as $h\to +\infty$) if for all $u \in X$ we have
\begin{itemize}
\item[{\rm (i)}] ({\it $\Gamma${-}liminf inequality}) for every
sequence $(u_h)_{h \in \N}$ converging to
$u$ in $X$,
$$
\liminf\limits _{h\to+\infty }F_h(u_h)\geq F(u);
$$
\item[{\rm (ii)}] ({\it $\Gamma${-}limsup inequality})
there exists a sequence
$(u_h)_{h \in \N}$ converging to $u$ in $X$, such that
$$
\limsup\limits _{h\to +\infty }F_h(u_h)\leq F(u).
$$
\end{itemize}
The function $F$ is called the $\Gamma${-}limit of
$(F_h)_{h \in \N}$ (with respect to $d$),
and we write $F\,=\,\Gamma{-}\lim_{h}F_h$.
$\Gamma${-}convergence is a convergence of
variational type as explained in the
following proposition.

\begin{proposition}
\label{Gamma-conv-prop}
Assume that the sequence $(F_h)_{h \in \N}$
$\Gamma${-}converges to $F$ and
that there exists a compact set $K\subseteq X$ 
such that for all $h \in \N$
$$
\inf\limits _{u\in K}F_h(u)=\inf\limits _{u\in X}
F_h(u).
$$
Then $F$ admits a minimum on $X$, $\inf_{X}F_h \to \min_X F$, and
any limit point of any sequence $(u_h)_{h \in \N}$ such that
$$
\lim\limits _{h\to +\infty }\Bigl( F_h(u_h)-
\inf\limits _{u\in X}F_h(u)\Bigr) =0,
$$
is a minimizer of $F$.
\end{proposition}

Moreover the following compactness result holds.

\begin{proposition}
\label{gconvcomp}
If $(X,d)$ is separable, and $(F_h)_{h \in \N}$ is a sequence of
functionals on $X$,
then there exists a subsequence $(F_{h_k})_{k \in \N}$ 
and a function
$F\,:\,X \to [-\infty;+\infty]$ such that $(F_{h_k})_{k \in \N}$
$\Gamma${-}converges to $F$.
\end{proposition}


\section{Discrete in time evolution of fractures in the cohesive case}
\label{discrevolsec}
In this section we are interested in generalized antiplanar shear of an
elastic body $\Om$ in
the context of linearized elasticity and in presence of cohesive fractures.
\par
The notion of {\it discrete in time evolution} for fractures relative to
time dependent boundary displacement
$g(t)$ and preexisting crack configuration $(\bar \Gamma,\bar \psi)$ has
been described in the Introduction.
It relies on the minimization of functionals of the form
\begin{equation}
\label{defprob1}
\|\nabla u\|^2+\int_{\Gamma \cup \Sg{g(t)}{u}}
\varphi(|[u]| \lor \psi)\,d\hs^{N-1},
\end{equation}
with $\psi$ positive function on $\Gamma$.
We now define rigorously the functional space to which the displacements
belong, and the properties
of $\Om$, $\Gamma$, $\psi$ and $g(t)$ in order to prove an existence result
for the discrete in time evolution
of fractures.
\par
Let $\Om$ be an open bounded subset of $\R^N$ with Lipschitz boundary.
Let $\partial_D \Om \subseteq \partial \Om$
be open in the relative topology, and let
$\partial_N \Om:=\partial \Om \setminus \partial_D \Om$.
Let $\varphi:[0,+\infty[ \to [0,+\infty[$ be increasing and concave,
$\varphi(0)=0$ and such that
$\lim_{s \to +\infty} \varphi(s)=1$. If $a:=\varphi'(0)<+\infty$, we have
\begin{equation}
\label{varphiabove2}
\varphi(s) \le as \quad\quad \text{for all }s \in [0,+\infty[.
\end{equation}
Let $T>0$, and let us consider a boundary displacement
$g \in AC([0,T];H^1(\Om))$
such that $\|g(t)\|_\infty \le C$ for all $t \in [0,T]$.
We discretize $g$ in the following way. Given $\delta>0$,
let $I_\delta$ be a subdivision
of $[0,T]$ of the form
$0=t^\delta_0<t^\delta_1<\dots<t^\delta_{N_\delta}=T$
such that $\max_i (t^\delta_i-t^\delta_{i-1})<\delta$.
For $0\le i \le N_\delta$ we set $g_i^\delta:=g(t_i^\delta)$.
\par
Let $\bar \Gamma \tsub \Om$ be rectifiable, 
and let $\bar \psi$ be a positive
function on
$\bar \Gamma$ such that
\begin{equation}
\label{crackfinite}
\int_{\bar \Gamma} \varphi(\bar \psi) \,d\hs^{N-1}<+\infty.
\end{equation}
Let us extend $\bar \psi$ to $\Omb$ setting 
$\bar \psi=0$ outside $\bar \Gamma$.
As for the space of the displacements, 
it would be natural following \cite{FL}
to consider $u \in SBV(\Om)$. Since $a=\varphi'(0)<+\infty$,
we have unfortunately that
the minimization of \eqref{defprob1} is not 
well posed in $SBV(\Om)$.
Let us in fact consider $(u_n)_{n \in \N}$ 
minimizing sequence for
\eqref{defprob1}: it turns out that we may 
assume $(u_n)_{n \in \N}$ bounded in $BV(\Om)$.
As a consequence $(u_n)_{n \in \N}$ admits 
a subsequence \wlystar convergent in
$BV(\Om)$ to a function $u \in BV(\Om)$. Then we have that
minimizing sequences of \eqref{defprob1}
converge (up to a subsequence) to a minimizer of 
the relaxation of
\eqref{defprob1} with
respect to the \wstar topology of $BV(\Om)$. By Proposition \ref{relaxation},
the natural domain of this relaxed functional is $BV(\Om)$,
and that its form is
\begin{equation}
\label{defen1rel}
\int_\Om f(\nabla u) \,dx+
\int_{\Gamma \cup \Sg{g(t)}{u}}
\varphi(|[u]| \lor \psi)\,d\hs^{N-1}+a|D^cu|(\Om),
\end{equation}
where
\begin{equation}
\label{bulkrelax2}
f(\xi):=
\begin{cases}
|\xi|^2 & \text{if } |\xi| \le \frac{a}{2}
\\ \\
\frac{a^2}{4}+a(|\xi|-\frac{a}{2}) & \text{if } |\xi| \ge \frac{a}{2}.
\end{cases}
\end{equation}
In view of these remarks, we consider $BV(\Om)$ as the space of
displacements $u$ of the body $\Om$, and a total energy of the form
\eqref{defen1rel}.
The volume part in the energy \eqref{defen1rel} can be interpreted
as the contribution of the elastic behavior of the body. The second
term represents the
work done to create the fracture $\Gamma \cup \Sg{g(t)}{u}$ with opening
given by $|[u]| \lor \psi$.
The new term $a|D^cu|$ can be interpreted as
the contribute of microcracks in the configuration which are considered
as reversible.
\par
Let us define the discrete evolution of the fracture in this new setting.
For $i=0$, let $u^\delta_0 \in BV(\Om)$ be a minimum of
\begin{equation}
\label{energyrelaxstep0}
\min_{u \in BV(\Om)} \left\{
\int_\Om f(\nabla u) \,dx+
\int_{\Sg{g^\delta_0}{u} \cup \bar \Gamma}
\varphi(|[u]|\lor \bar \psi)\,d\hs^{N-1}
+a|D^cu|
\right\}.
\end{equation}
We set $\Gamma^\delta_0:=\Sg{g^\delta_0}{u^\delta_0} \cup \bar \Gamma$.
\par
Supposing to have constructed $u^\delta_j$ and $\Gamma^\delta_j$ for
all $j=0,\dots,i-1$, let
$u^\delta_i$ be a minimum of
\begin{equation}
\label{energyrelaxstepi}
\min_{u \in BV(\Om)} \left\{
\int_\Om f(\nabla u) \,dx+
\int_{\Sg{g^\delta_i}{u} \cup \Gamma^\delta_{i-1}}
\varphi(|[u]| \lor \psi^\delta_{i-1})\,d\hs^{N-1}
+a|D^cu|,
\right\}
\end{equation}
where $\psi^\delta_{i-1}:=\bar \psi \lor |[u^\delta_0]|
\lor \dots \lor |[u^\delta_{i-1}]|$. We set
$\Gamma^\delta_i:=\Gamma^\delta_{i-1} \cup \Sg{g^\delta_i}{u^\delta_i}$.
\par
The following proposition establish the existence of this discrete evolution.

\begin{proposition}
\label{discrevol}
Let $I_\delta=\{0=t^\delta_0<\dots<t^\delta_{N_\delta}=T\}$ be a 
subdivision of $[0,T]$
such that $\max(t^\delta_i-t^\delta_{i-1})<\delta$, let 
$\bar \Gamma$ be a preexisting crack,
and $\bar \psi$ a positive function on $\bar \Gamma$ satisfying 
\eqref{crackfinite} and
extended to zero outside $\bar \Gamma$.
Then for all $i=0,\dots,N_\delta$ there exists $u^\delta_i \in BV(\Om)$ 
such that setting
$\Gamma^\delta_{-1}:=\bar \Gamma$, $\psi^\delta_{-1}:=\bar \psi$ and
\begin{equation}
\label{discrfracture}
\Gamma^\delta_i:= \bar \Gamma \cup \bigcup_{j=0}^{i} \Sg{g^\delta_j}{u^\delta_j},
\quad\quad
\psi^\delta_i(x):=\bar \psi(x) \lor |[u^\delta_0]|(x) \lor \dots 
\lor |[u^\delta_i]|(x)
\end{equation}
the following holds:
\begin{itemize}
\item[(a)] $\|u^\delta_i\|_\infty \le \|g^\delta_i\|_\infty \le C$;
\item[]
\item[(b)] for all $v \in BV(\Om)$ we have
\begin{multline}
\label{comparison1}
\int_\Om f(\nabla u^\delta_i) \,dx+
\int_{\Gamma^\delta_i} \varphi(\psi^\delta_i)\,d\hs^{N-1}
+a|D^cu^\delta_i|(\Om) \\
\le
\int_\Om f(\nabla v) \,dx+
\int_{\Sg{g^\delta_i}{v} \cup
\Gamma^\delta_{i-1}} \varphi(|[v]| \lor \psi^\delta_{i-1})\,d\hs^{N-1} 
+a|D^cv|(\Om),
\end{multline}
where $a=\varphi'(0)$ and $f$ is defined in \eqref{bulkrelax2};
\item[]
\item[(c)] we have that
\begin{multline}
\label{sbvenergy}
\int_\Om f(\nabla u^\delta_i) \,dx+
\int_{\Gamma^\delta_i} \varphi(\psi^\delta_i)\,d\hs^{N-1}
+a|D^cu^\delta_i|(\Om) \\
=\inf_{v \in SBV(\Om)}
\left\{
\|\nabla v\|^2+\int_{\Sg{g^\delta_i}{v} \cup
\Gamma^\delta_{i-1}} \varphi(|[v]| \lor \psi^\delta_{i-1})\,d\hs^{N-1}
\right\}.
\end{multline}
\end{itemize}
\end{proposition}

\begin{proof}
We have to prove that problems \eqref{energyrelaxstep0} and 
\eqref{energyrelaxstepi} admit solutions.
Let us consider for example problem \eqref{energyrelaxstepi}, 
the other being similar.
Let $(u_n)_{n \in \N}$ be a minimizing sequence for problem 
\eqref{energyrelaxstepi}.
By a truncation argument we may assume that 
$\|u_n\|_\infty \le \|g^\delta_i\|$.
Comparing $u_n$ with $g^\delta_i$, we get for $n$ large
\begin{multline}
\int_\Om f(\nabla u_n) \,dx+
\int_{\Sg{g^\delta_i}{u_n} \cup \Gamma^\delta_{i-1}} \varphi(|[u_n]|
\lor \psi^\delta_{i-1})\,d\hs^{N-1} +a|D^cu_n|(\Om) \\
\le \int_\Om f(\nabla g^\delta_0)\,dx
+\int_{\Gamma^\delta_{i-1}}
\varphi(\psi^\delta_{i-1})\,d\hs^{N-1}+1 \le C',
\end{multline}
with $C'$ independent of $n$. Since there exists $d>0$ such that
$a|\xi|-d \le f(\xi)$ for all $\xi \in \R^N$,
we deduce that $(\nabla u_n)_{n \in \N}$ 
is bounded in $L^1(\Om;\R^N)$.
Moreover if $\bar s$ is such that $\varphi(\bar s)=\frac{1}{2}$ 
and $\bar a$ is such that $s \le \bar a \varphi(s)$
for all $s \in [0, \bar s]$, we have
\begin{align}
\int_{S(u_n)}|[u_n]|\,d\hs^{N-1} &=
\int_{\{|[u_n]| \le \bar s\}}|[u_n]|\,d\hs^{N-1} 
+\|g^\delta_i\|_\infty \hs^{N-1}(\{|[u_n]|>\bar s\}) \\
\nonumber
&\le \bar a \int_{|[u_n]| \le \bar s} \varphi(|[u_n]|)\,d\hs^{N-1}
+2\|g^\delta_i\|_\infty
\int_{|[u_n]|>\bar s} \varphi(|[u_n]|)\,d\hs^{N-1} \\
\nonumber
&\le (\bar a+2\|g^\delta_i\|_\infty)C'.
\end{align}
Finally for all $n$
$$
|D^cu_n| \le \frac{C'}{a}.
$$
We conclude that $(u_n)_{n \in \N}$ 
is bounded in $BV(\Om)$. Then there exists 
$u \in BV(\Om)$ such that
up to a subsequence $u_n \weakst u$ \wlystar in $BV(\Om)$ and 
pointwise almost everywhere.
Let us set $u^\delta_i:=u$. By Lemma \ref{lsc} we deduce that
\begin{multline}
\int_\Om f(\nabla u) \,dx+
\int_{\Sg{g^\delta-i}{u} \cup \Gamma^\delta_{i-1}} \varphi(|[u]|
\lor \psi^\delta_{i-1})\,d\hs^{N-1} +a|D^cu|(\Om) \\
\le \liminf_n
\int_\Om f(\nabla u_n) \,dx+
\int_{\Sg{g^\delta_i}{u_n} \cup \Gamma^\delta_{i-1}} \varphi(|[u_n]|
\lor \psi^\delta_{i-1})\,d\hs^{N-1} +a|D^cu_n|(\Om).
\end{multline}
Setting $\psi^\delta_i:=\psi^\delta_{i-1} \lor |[u^\delta_i]|$,
we have that point $(b)$ holds. Moreover $\|u^\delta_i\|_\infty \le 
\|g^\delta_0\|_\infty \le C$,
so that point $(a)$ holds. Finally point $(c)$ is a consequence of 
Proposition \ref{relaxation}.
\end{proof}

Let us consider now the following piecewise constant interpolation 
in time:
\begin{equation}
\label{interpbis}
u^\delta(t):=u^\delta_i, \quad\quad
\Gamma^\delta(t):=\Gamma^\delta_i, \quad\quad
\psi^\delta(t):=\psi^\delta_i, \quad\quad g^\delta(t):=g^\delta_i
\quad\quad
t^\delta_i \le t <t^\delta_{i+1}
\end{equation}
with $u^\delta(T):=u^\delta_{N_\delta}$, 
$\Gamma^\delta(T):=\Gamma^\delta_{N_\delta}$,
$\psi^\delta(T):=\psi^\delta_{N_\delta}$, and $g^\delta(T):=g(T)$.
\par
For every $v \in BV(\Om)$ and for every $t \in [0,T]$ let us set
\begin{equation}
\label{energydefi}
\eub^\delta(t,v):=
\int_\Om f(\nabla v) \,dx+
\int_{\Sg{g^\delta(t)}{v} \cup
\Gamma^\delta(t)} \varphi(|[v]| \lor \psi^\delta(t))\,d\hs^{N-1}
+a|D^cv|(\Om).
\end{equation}
Then the following estimate holds.

\begin{lemma}
\label{energyabove}
There exists $e^\delta_a \to 0$ for $\delta \to 0$ and $a \to +\infty$
such that for all $t \in [0,T]$ we have
\begin{equation}
\label{enabove}
\eub^\delta(t,u^\delta(t)) \le
\eub^\delta(0,u^\delta(0))+
\int_{0}^{t^\delta_i} \int_\Om f'(\nabla u^\delta(\tau))
\nabla \dot{g}(\tau)\,dx\,d\tau+e^\delta_a,
\end{equation}
where $t^{\delta}_{i}$ is the step discretization point such that
$t^{\delta}_{i} \le t<t^{\delta}_{i+1}$.
\end{lemma}

\begin{proof}
Comparing $u^\delta_i$ with $u^\delta_{i-1}+g^\delta_i-g^\delta_{i-1}$ 
by means of \eqref{comparison1} we
obtain
$$
\eub^\delta(t^\delta_i,u^\delta_i)
\le \int_{\Om} f(\nabla u^\delta_{i-1}
+\nabla g^\delta_i-\nabla g^\delta_{i-1}) \,dx
+\int_{\Gamma^\delta_{i-1}} \varphi(\psi^\delta_{i-1})\,d\hs^{N-1} 
+a|D^c u^\delta_{i-1}|(\Om).
$$
Notice that by the very definition of $f$ the following hold:
\begin{itemize}
\item[1)]
if $|\nabla u^\delta_{i-1}+\nabla g^\delta_i-\nabla g^\delta_{i-1}| 
\ge \frac{a}{2}$ and
$|\nabla u^\delta_{i-1}| \ge \frac{a}{2}$
$$
f'(\nabla u^\delta_{i-1}+\nabla g^\delta_i-\nabla g^\delta_{i-1})
=f'(\nabla u^\delta_{i-1});
$$
\item[]
\item[2)] if $|\nabla u^\delta_{i-1}+
\nabla g^\delta_i-\nabla g^\delta_{i-1}| < \frac{a}{2}$ and
$|\nabla u^\delta_{i-1}| \ge \frac{a}{2}$
$$
f(\nabla u^\delta_{i-1}+\nabla g^\delta_i-\nabla g^\delta_{i-1}) 
\le f(\nabla u^\delta_{i-1});
$$
\item[]
\item[3)]
if $|\nabla u^\delta_{i-1}+\nabla g^\delta_i-\nabla g^\delta_{i-1}|
\ge \frac{a}{2}$ and
$|\nabla u^\delta_{i-1}| < \frac{a}{2}$
$$
f(\nabla u^\delta_{i-1}+\nabla g^\delta_i-\nabla g^\delta_{i-1}) \le
f(\nabla u^\delta_{i-1})
+2(\nabla u^\delta_{i-1},\nabla g^\delta_i-\nabla g^\delta_{i-1})+
|\nabla g^\delta_i-\nabla g^\delta_{i-1}|^2;
$$
\item[]
\item[4)]
if $|\nabla u^\delta_{i-1}+\nabla g^\delta_i-\nabla g^\delta_{i-1}|
< \frac{a}{2}$ and
$|\nabla u^\delta_{i-1}| < \frac{a}{2}$
$$
f(\nabla u^\delta_{i-1}+\nabla g^\delta_i-\nabla g^\delta_{i-1})=
f(\nabla u^\delta_{i-1})
+2(\nabla u^\delta_{i-1},\nabla g^\delta_i-\nabla g^\delta_{i-1})+
|\nabla g^\delta_i-\nabla g^\delta_{i-1}|^2.
$$
\end{itemize}
Then by convexity of $f$ we deduce
$$
\eub^\delta(t^\delta_i,u^\delta_i)
\le
\eub^\delta(t^\delta_{i-1},u^\delta_{i-1})
+\int_{\Om} f'(\nabla u^\delta_{i-1})
(\nabla g^\delta_i-\nabla g^\delta_{i-1})\,dx+
R^{\delta,a}_{i-1},
$$
where
$$
R^{\delta,a}_{i-1}:= \int_{\Om}
|\nabla g^\delta_i-\nabla g^\delta_{i-1}|^2\,dx+
\int_{\{|\nabla u^\delta_{i-1}| \ge \frac{a}{2}\}}
|f'(\nabla u^\delta_{i-1})|\,
|\nabla g^\delta_i-\nabla g^\delta_{i-1}|\,dx.
$$
Then summing up from $t^\delta_i$ to $t^\delta_0$, and taking into account
\eqref{interpbis} we get
$$
\eub^\delta(t,u^\delta(t))
\le
\eub^\delta(0,u^\delta(0))
+\int_{0}^{t^\delta_i} \int_{\Om} f'(\nabla u^\delta(\tau))
\nabla \dot{g}(\tau)\,dx\,d\tau+ \int_{0}^{t^\delta_i} 
R^{\delta,a}(\tau)\,d\tau,
$$
where
\begin{equation}
\label{errorR}
R^{\delta,a}(\tau):= \sigma(\delta)\|\nabla \dot{g}(\tau)\|+
\int_{\{|\nabla u^\delta(\tau)| \ge \frac{a}{2}\}} 
|f'(\nabla u^\delta(\tau))|\,
|\nabla \dot{g}(\tau)|\,dx
\end{equation}
and
$$
\sigma(\delta):=\max_{i=1,\dots,N_\delta} 
\int_{t^\delta_{i-1}}^{t^\delta_i}\|\nabla \dot{g}\|\,d\tau.
$$
In order to conclude the proof it is sufficient to see that
$$
\int_0^T R^{\delta,a}(\tau)\,d\tau \to 0
$$
as $\delta \to 0$ and $a \to +\infty$. Notice that 
$\sigma(\delta) \to 0$ as $\delta \to 0$
by the absolutely continuity of $\|\nabla \dot{g}\|$. 
Let us come to the second term.
Notice that $|f'(\nabla u^\delta(\tau))|=a$ on
$\{|\nabla u^\delta(\tau)| \ge \frac{a}{2}\}$.
Then we have to see
\begin{equation}
\label{Rto0}
\int_0^T \int_\Om a|\nabla \dot{g}(\tau)|
\1_{\{|\nabla u^\delta(\tau)| \ge \frac{a}{2}\}}\,dx\,d\tau \to 0
\end{equation}
as $\delta \to 0$ and $a\to +\infty$. Setting
$A^\delta_a(\tau):=\{x \in \Om\,:\,|\nabla u^\delta(\tau)|(x) \ge \frac{a}{2}\}$
we have by H\"older inequality
$$
\int_\Om a|\nabla \dot{g}(\tau)| \1_{A^\delta_a(\tau)}\,dx \le
a\sqrt{|A^\delta_a(\tau)|} \left( \int_{A^\delta_a(\tau)}
|\nabla \dot{g}(\tau)|^2\,dx \right)^{\frac{1}{2}}.
$$
Notice that
\begin{equation}
\label{estadeltac}
\frac{a^2}{2}|A^\delta_a(\tau)| \le a\int_{A^\delta_a(\tau)}
|\nabla u^\delta(\tau)|\,dx \le
2\int_{A^\delta_a(\tau)} f(\nabla u^\delta(\tau))\,dx \le C',
\end{equation}
where $C'$ depends only on $g$ and is obtained comparing
$u^\delta(\tau)$ with $g^\delta(\tau)$ by means of
\eqref{comparison1}.
We deduce that
\begin{equation}
\label{finest}
\int_\Om a|\nabla \dot{g}(\tau)| \1_{A^\delta_a(\tau)}\,dx \le
\sqrt{2C'}\left( \int_{A^\delta_a(\tau)}|\nabla \dot{g}(\tau)|^2\,dx
\right)^{\frac{1}{2}} \le
\sqrt{2C'} \|\nabla \dot{g}(\tau)\|.
\end{equation}
As $\delta \to 0$ and $a \to +\infty$, by \eqref{estadeltac}
we have that $|A^\delta_a(\tau)| \to 0$.
Then by the equicontinuity of $|\nabla \dot{g}(\tau)|^2$ and
by the Dominated Convergence Theorem,
we deduce that \eqref{Rto0} holds, and the proof is finished.
\end{proof}

\section{The main results}
\label{mainressec}
Let $\Om$ be an open bounded subset of $\R^N$ with Lipschitz 
boundary.
Let $\partial_D \Om \subseteq \partial \Om$ be open in the 
relative topology,
and let $\partial_N \Om:=\partial \Om \setminus \partial_D \Om$.
\par
In this section we consider discrete in time evolution of fractures
in a linearly elastic body whose reference configuration is 
given by $\Om_h:=h\Om$,
where $h>0$. Let us assume that the cohesive forces on the
fractures of $\Om_h$ are given in the sense of Section 
\ref{discrevolsec}
by a function $\varphi\,:\,[0,+\infty[ \to [0,1]$
which is increasing, concave,
$\varphi(0)=0$, $\varphi'(0)=a<+\infty$ and such that
$\lim_{s \to +\infty} \varphi(s)=1$. Let us
moreover set
\begin{equation}
\label{bulkrelax3}
f(\xi):=
\begin{cases}
|\xi|^2 & \text{if } |\xi| \le \frac{a}{2}
\\ \\
\frac{a^2}{4}+a(|\xi|-\frac{a}{2}) & \text{if }
|\xi| \ge \frac{a}{2}.
\end{cases}
\end{equation}
Let us consider on $\partial_D \Om_h:=h \partial_D \Om$
boundary displacements
of the following particular form
\begin{equation}
\label{bdrydisp}
g_h(t,x):=h^\alpha g\left( t,\frac{x}{h} \right)
\end{equation}
with $g \in AC([0,T];H^1(\Om))$ such that
$\|g(t)\|_\infty \le C$ for all $t \in [0,T]$.
Let moreover $\bar \Gamma \tsub \Om$
be rectifiable
with 
\begin{equation}
\label{gbarfinite}
\hs^{N-1}(\bar \Gamma)<+\infty
\end{equation}
and let $\bar \gamma$
be a positive function defined on $\bar \Gamma$. 
We extend $\bar \gamma$ to $\Omb$ setting $\bar \gamma=0$ outside
$\bar \Gamma$. Let us consider
$({\bar \Gamma}_h,\bar \psi_h)$ as a
preexisting crack configuration in $\Om_h$, where
\begin{equation}
\label{initconfig}
{\bar \Gamma}_h:=h\bar \Gamma,
\quad\quad
\bar \psi_h(x):=\bar \gamma \left( \frac{x}{h} \right),
\quad x \in \Om_h.
\end{equation}
Given $\delta>0$, let 
$I_\delta=\{0=t^\delta_0<\dots<t^\delta_{N_\delta}=T\}$ be a
subdivision of $[0,T]$ such that
$\max (t^\delta_i-t^\delta_{i-1})<\delta$, and let
$\{t \to (u^{\delta,h}(t), \Gamma^{\delta,h}(t),\psi^{\delta,h}(t))
\,:\,t \in [0,T]\}$
be the piecewise constant interpolation
in the sense of \eqref{interpbis} of a discrete in time
evolution of fractures in $\Om_h$
relative to the boundary datum $g_h$, the preexisting crack
configuration $({\bar \Gamma}_h, {\bar \psi}_h)$
and the subdivision $I_\delta$ given by Proposition \ref{discrevol}.
\par
Our aim is to study the asymptotic behavior of
$\{t \to (u^{\delta,h}(t),\Gamma^{\delta,h}(t),\psi^{\delta,h}(t))\,:\,
t \in [0,T]\}$ as
$\delta \to 0$ and $h \to +\infty$. Let us
consider $h \in \N$ (we can consider any sequence 
which diverges to $+\infty$), let us fix $\delta_h \to 0$,
and let us set for all $t \in [0,T]$
\begin{equation}
\label{defuh}
u_h(t):=u^{\delta_h,h}(t),
\quad\quad
\Gamma_h(t):=\Gamma^{\delta_h,h}(t),
\quad\quad
\psi_h(t):=\psi^{\delta_h,h}(t),
\end{equation}
and let $g^\delta_h(t):=g_h(t^\delta_i)$ where 
$t^\delta_i \in I_\delta$ is such that $t^\delta_i
\le t < t^\delta_{i+1}$.
Let us moreover set for every $v \in BV(\Om)$ and 
for every $t \in [0,T]$
\begin{equation}
\label{energydefth}
\eub_h(t,v):=
\int_\Om f(\nabla v) \,dx+
\int_{\Sg{g^\delta_h(t)}{v} \cup
\Gamma_h(t)} \varphi(|[v]| \lor \psi_h(t))\,d\hs^{N-1} 
+a|D^cv|(\Om).
\end{equation}
The asymptotic of $(u_h,\Gamma_h,\psi_h)$ depends on $\alpha$,
and we have to distinguish three cases.
The first case $\alpha=\frac{1}{2}$ was stated in the 
Introduction and reveals the prevalence
of brittle effects as the size of the body increases.
We give here the precise statement
we will prove.

\begin{theorem}
\label{main1}
Let $g \in AC(0,T;H^1(\Om))$ be such that $\|g(t)\|_\infty \le C$
for all $t \in [0,T]$.
Let $\{t \to (u_h(t),\Gamma_h(t),\psi_h(t))\,:\,t \in [0,T]\}$
be the piecewise constant interpolation of a discrete in time
evolution of fractures in $\Om_h$ relative to the
preexisting crack configuration $({\bar \Gamma}_h,\bar \psi_h)$ 
and the boundary data
$$
g_h(x,t):=\sqrt{h}g \left( \frac{x}{h},t \right).
$$
Then the following facts hold:
\begin{itemize}
\item[]
\item[(a)] there exists a constant $C'$ dependent only on $g$ 
such that for all $t \in [0,T]$
\begin{equation}
\label{energybound1}
\frac{1}{h^{N-1}} \eub_h(t,u_h(t)) \le C';
\end{equation}
\item[]
\item[(b)] for all $t \in [0,T]$
\begin{equation}
\label{dispconv1}
v_h(t,x):=\frac{1}{\sqrt{h}} u_h(t,hx) \quad\quad 
\text{is bounded in }BV(\Om);
\end{equation}
\item[(c)] there exists a subsequence independent of $t$ and
there exists a quasistatic evolution of brittle fractures
$\{t \to (v(t), K(t))\,:\,t \in [0,T]\}$
in $\Om$ relative to the preexisting crack $\bar \Gamma$ 
and boundary displacement $g$
in the sense of Theorem \ref{qse} such that for all $t \in [0,T]$
we have
\begin{equation}
\label{gradconv1}
\nabla v_h(t) \weak \nabla v(t) \quad\quad
\text{weakly in }L^1(\Om;\R^N),
\end{equation}
and every accumulation point $v$ of $(v_h(t))_{h \in \N}$
in the \wstar topology of $BV(\Om)$ is such that
$v \in SBV(\Om)$, $\Sg{g(t)}{v} \tsub K(t)$ and
$\nabla v=\nabla v(t)$. Moreover
for all $t \in [0,T]$ we have
\begin{equation}
\label{totconv1}
\frac{1}{h^{N-1}} \eub_h(t,u_h(t)) \to
\|\nabla v(t)\|^2+\hs^{N-1}(K(t));
\end{equation}
in particular $h^{-N+1}|D^cu_h(t)|(\Om_h)| \to 0$,
\begin{equation}
\label{bulkconv1}
\frac{1}{h^{N-1}}\int_{\Om_h} f(\nabla u_h(t))\,dx \to 
\|\nabla v(t)\|^2,
\end{equation}
and
\begin{equation}
\label{fracconv1}
\frac{1}{h^{N-1}}\int_{\Gamma_h(t)} \varphi(\psi_h(t)) \,d\hs^{N-1} 
\to \hs^{N-1}(K(t)).
\end{equation}
\end{itemize}
\end{theorem}

The case $\alpha<\frac{1}{2}$ leads to a problem in elasticity 
in $\Om_h \setminus {\bar \Gamma}_h$ in the sense of 
the following theorem.

\begin{theorem}
\label{main2}
Let $g \in AC(0,T;H^1(\Om))$ be such that $\|g(t)\|_\infty \le C$
for all $t \in [0,T]$.
Let $\{t \to (u_h(t),\Gamma_h(t),\psi_h(t))\,:\,t \in [0,T]\}$
be the piecewise constant interpolation of a discrete in time
evolution of fractures in $\Om_h$ relative to the
initial crack configuration $({\bar \Gamma}_h,\bar \psi_h)$ and
the boundary data
$$
g_h(x,t):=h^{\alpha}g \left( t,\frac{x}{h} \right)
$$
with $\alpha<\frac{1}{2}$.
Then the following facts hold:
\begin{itemize}
\item[]
\item[(a)] for all $t \in [0,T]$
\begin{equation}
\label{dispconv2}
v_h(t,x):=\frac{1}{h^\alpha} u_h(t,hx) \quad\quad
\text{is bounded in }BV(\Om);
\end{equation}
\item[]
\item[(b)] for all $t \in [0,T]$ every accumulation point $v(t)$
of $(v_h(t))_{h \in\N}$ in
the \wstar topology of $BV(\Om)$ is such that $v(t) \in SBV(\Om)$ and
$\Sg{g(t)}{v(t)} \tsub \bar \Gamma$;
\item[]
\item[(c)]
if $\varphi(s)=1$ for $s \ge \bar s$, and $\bar \gamma \ge \eps>0$,
then there exists
a subsequence independent of $t$ such that for all $t \in [0,T]$
\begin{equation}
\label{gradconv2}
\nabla v_h(t,x) \weak \nabla v(t) \quad\quad
\text{weakly in }L^1(\Om;\R^N),
\end{equation}
where $v(t)$ is a minimizer of
\begin{equation}
\label{elprob}
\min \{\|\nabla v\|^2\,:\, v \in SBV(\Om),\,\Sg{g(t)}{v}
\tsub \bar \Gamma\};
\end{equation}
moreover for all $t \in [0,T]$ we have
\begin{equation}
\label{totconv2}
\frac{1}{h^{N+2\alpha-2}} \int_{\Om_h}f(\nabla u_h(t))\,dx \to
\|\nabla v(t)\|^2.
\end{equation}
\end{itemize}
\end{theorem}

Finally for the case $\alpha>\frac{1}{2}$ the body goes to
{\it rupture} at time $t=0$, in the
sense of the following theorem.

\begin{theorem}
\label{main3}
Let $g \in AC(0,T;H^1(\Om))$ be such that $\|g(t)\|_\infty \le C$
for all $t \in [0,T]$.
Let $\{t \to (u_h(t),\Gamma_h(t),\psi_h(t))\,:\,t \in [0,T]\}$
be the piecewise constant interpolation of a discrete in time
evolution of fractures in $\Om_h$ relative to the
initial crack configuration $({\bar \Gamma}_h,\bar \psi_h)$ and
the boundary data
$$
g_h(x,t):=h^{\alpha}g \left( \frac{x}{h},t \right)
$$
with $\alpha>\frac{1}{2}$. Let us set $v_h(t,x):=\frac{1}{h^\alpha}
u_h(t,hx)$ for all $x \in \Om$ and
for all $t \in [0,T]$.
\par
Then $(v_h(0))_{h \in \N}$
is bounded in $BV(\Om)$, and every accumulation
point $v$ of
$(v_h(0))_{h \in \N}$ in the \wstar 
topology of $BV(\Om)$ is piecewise constant in $\Om$, that is
$v \in SBV(\Om)$ and $\nabla v=0$. Moreover
\begin{equation}
\label{minpiece}
\hs^{N-1}(\Sg{g(0)}{v(0)} \cup \bar \Gamma) \le 
\hs^{N-1}(\Sg{g(0)}{w} \cup \bar \Gamma)
\end{equation}
for all piecewise constant function $w \in SBV(\Om)$.
\end{theorem}

\section{Proof of Theorem \ref{main1}}
\label{main1sec}
In this section we will give the proof of Theorem \ref{main1}.
Let $\{t \to (u_h(t),\Gamma_h(t),\psi_h(t))\,:\,t \in [0,T]\}$
be the piecewise constant interpolation
of a discrete in time evolution of cohesive fracture in $\Om_h$
relative the subdivision
$I_{\delta_h}:=\{0=t^{\delta_h}_0<\dots<t^{\delta_h}_{N_{\delta_h}}=T\}$, 
the preexisting crack configuration
$({\bar \Gamma}_h,{\bar \psi}_h)$ given by \eqref{initconfig}
and the boundary displacement $\sqrt{t}g(t,\frac{x}{h})$.
\par
In order to prove Theorem \ref{main1}, we need some preliminary
analysis.
First of all, it is convenient to rescale $u_h$ and $\Gamma_h$ in
the following way: for all
$t \in [0,T]$ let $v_h(t) \in BV(\Om)$ and $K_h(t) \tsub \Om
\cup \partial_D \Om$ be defined by
\begin{equation}
\label{riscaledisp1}
v_h(t,x):=\frac{1}{\sqrt{h}}u_h \left( t,hx \right),
\quad\quad
K_h(t):=\frac{1}{h} \Gamma_h(t).
\end{equation}
Let us moreover set
\begin{equation}
\label{defgammapiccoloh}
\gamma_h(t,x):=\frac{1}{\sqrt{h}} \psi_h(t,hx)=\max_{s \le t}
|[v_h(s)](t,x)| \lor \bar \gamma(x),
\quad t \in [0,T], x \in \Om.
\end{equation}
We notice that $\{t \to (v_h(t),K_h(t),\gamma_h(t))\,:\,t \in [0,T]\}$
is the piecewise constant interpolation
of a discrete in time evolution of cohesive fractures in $\Om$ relative
to the subdivision $I_{\delta_h}$,
the preexisting crack configuration $(\bar \Gamma,\bar \gamma)$ and
boundary displacement $g(t)$ with
respect to the basic total energy
\begin{equation}
\label{riscenergy1}
\int_\Om f_h(\nabla v) \,dx+
\int_{\Sg{g^{\delta_h}(t)}{v} \cup K_h(t)} \varphi_{h}(|[v]|
\lor \gamma_h(t))\,d\hs^{N-1}+a \sqrt{h}|D^cv|(\Om),
\end{equation}
where
\begin{equation}
\label{defvarphi1h}
\varphi_h(s):=\varphi(\sqrt{h}s),
\end{equation}
and
\begin{equation}
\label{deff1h}
f_h(\xi):=
\begin{cases}
|\xi|^2 & \text{if } |\xi| \le \frac{a \sqrt{h}}{2}
\\ \\
\frac{a^2h}{4}+a \sqrt{h}(|\xi|-\frac{a \sqrt{h}}{2}) & \text{if }
|\xi| \ge \frac{a \sqrt{h}}{2}.
\end{cases}
\end{equation}
\vskip10pt
Let us recall some properties of the evolution
$\{t \to (v_h(t),K_h(t),\gamma_h(t))\,:\,t \in [0,T]\}$
which derive from Proposition \ref{discrevol} and
that will be employed in the sequel:
\begin{itemize}
\item[(a)] for all $t \in [0,T]$
\begin{equation}
\label{infty1}
\|v_h(t)\|_\infty \le \|g^{\delta_h}(t)\|_\infty;
\end{equation}
\item[]
\item[(b)] for all $w \in BV(\Om)$ we have
\begin{multline}
\label{comparezero1}
\int_{\Om} f_h(\nabla v_h(0)) \,dx
+\int_{K_h(0)} \varphi_h(|[v_h(0)]| \lor \bar \gamma)\,d\hs^{N-1}
+a\sqrt{h}|D^c v_h(0)|(\Om) \\
\le \int_{\Om} f_h(\nabla w) \,dx
+\int_{\Sg{g^{\delta_h}(0)}{w} \cup \bar \Gamma}
\varphi_h(|[w]| \lor \bar \gamma)\,d\hs^{N-1} +a \sqrt{h}|D^cw|(\Om);
\end{multline}
\item[]
\item[(c)] for all $w \in BV(\Om)$ and for all $t \in ]0,T]$ we have
\begin{multline}
\label{compare1}
\int_{\Om} f_h(\nabla v_h(t)) \,dx
+\int_{K_h(t)} \varphi_h(\gamma_h(t))\,d\hs^{N-1}
+a\sqrt{h}|D^c v_h(t)|(\Om) \\
\le \int_{\Om} f_h(\nabla w) \,dx
+\int_{\Sg{g^{\delta_h}(t)}{w} \cup K_h(t)}
\varphi_h(|[w]| \lor \gamma_h(t))\,d\hs^{N-1} +a \sqrt{h}|D^cw|(\Om).
\end{multline}
\end{itemize}
Let us set for all $v \in BV(\Om)$ and for all $t \in [0,T]$
\begin{equation}
\label{deffs}
\fs_h(t,w):=\int_{\Om} f_h(\nabla w) \,dx
+\int_{\Sg{g^{\delta_h}(t)}{w} \cup K_h(t)}
\varphi_h(|[w]| \lor \gamma_h(t))\,d\hs^{N-1} +a \sqrt{h}|D^cw|(\Om).
\end{equation}
Notice that for all $t \in [0,T]$
\begin{equation}
\label{fseub}
\fs_h(t,v_h(t))=\frac{1}{h^{N-1}} \eub_h(t,u_h(t)),
\end{equation}
where $\eub_h(t,u)$ is defined in \eqref{energydefth}.
\par
Recalling Lemma \ref{energyabove}, we have that the following holds.

\begin{lemma}
\label{energy1above}
For all $t \in [0,T]$ we have
\begin{equation}
\label{en1above}
\fs_h(t,v_h(t)) \le \fs_h(0,v_h(0))+
\int_{0}^{t_h} \int_\Om f'_h(\nabla v_h(\tau)) \nabla
\dot{g}(\tau)\,dx\,d\tau+e(h),
\end{equation}
where $e(h) \to 0$ as $h \to +\infty$, and $t_h:=t^{\delta_h}_{i_h}$
is the step
discretization point of $I_{\delta_h}$ such that
$t^{\delta_h}_{i_h} \le t<t^{\delta_h}_{i_h+1}$.
\end{lemma}

The following corollary provides a bound on the total energy of the
discrete in time evolution.

\begin{corollary}
\label{enboundh1lem}
There exists a constant $C'$ independent of $h$ such that for all
$t \in [0,T]$ we have
\begin{equation}
\label{ebound1h}
\fs_h(t,v_h(t))+\|v_h(t)\|_\infty \le C'.
\end{equation}
\end{corollary}

\begin{proof}
By \eqref{en1above} we have
$$
\fs_h(t,v_h(t)) \le \fs_h(0,v_h(0))
+\int_{0}^{t_h} \int_\Om f'_h(\nabla v_h(\tau))
\nabla \dot{g}(\tau)\,dx\,d\tau+e(h),
$$
where $e(h) \to 0$ as $h \to +\infty$, and $t_h:=t^{\delta_h}_{i_h}$
is such that
$t^{\delta_h}_{i_h} \le t<t^{\delta_h}_{i_h+1}$.
\par
Notice that by \eqref{comparezero1} we have
$$
\fs_h(0,u_h(0)) \le \|\nabla g(0)\|^2+\hs^{N-1}(\bar \Gamma),
$$
and by \eqref{compare1} for all $\tau \in [0,T]$
$$
\int_\Om f_h(\nabla v_h(\tau))\,dx \le \|\nabla g^{\delta_h}(\tau)\|^2.
$$
Moreover for all $\tau \in [0,T]$
$$
\int_\Om |f'_h(\nabla v_h(\tau))|^2 \,dx \le 4
\int_\Om f_h(\nabla v_h(\tau))\,dx
$$
and so, taking into account \eqref{infty1},
we deduce that \eqref{ebound1h} holds.
\end{proof}

As a consequence of Corollary \ref{enboundh1lem}, we infer a uniform
bound on the total variations of $v_h(t)$.

\begin{corollary}
\label{totvarbound1h}
There exists $C''$ independent of
$h$ such that for all $t \in [0,T]$ we have
\begin{equation}
\label{varbound}
|Dv_h(t)|(\Om) \le C''.
\end{equation}
\end{corollary}

\begin{proof}
In fact, since for $h$ large we have for all $\xi \in \R^N$
$$
|\xi|-1 \le f_h(\xi),
$$
we deduce that for all $t \in [0,T]$
$$
\int_\Om |\nabla v_h(t)|\,dx \le
\int_\Om [f_h(\nabla v_h(t))+1]\,dx \le C'+|\Om|,
$$
where $C'$ is given by Corollary \ref{enboundh1lem}.
Moreover if $\bar s$ is such that $\varphi(\bar s)=\frac{1}{2}$ and
$\bar a$ is such that $s \le \bar a \varphi(s)$
for all $s \in [0, \bar s]$, we have for all $h$ and 
for all $t \in [0,T]$
\begin{multline}
\int_{S(v_h(t))}|[v_h(t)]|\,d\hs^{N-1} =
\int_{|[v_h(t)]| \le \frac{\bar s}{\sqrt{h}}}|[v_h(t)]|\,d\hs^{N-1}
+\|v_h(t)\|_\infty
\hs^{N-1} \left( \Big\{ |[v_h(t)]| \ge
\frac{\bar s}{\sqrt{h}} \Big\} \right) \\
\le \bar a \int_{|[v_h(t)]| \le \frac{\bar s}{\sqrt{h}}}
\varphi_h(|[v_h(t)]|)\,d\hs^{N-1}+2C'
\int_{|[v_h(t)]| \ge \frac{\bar s}{\sqrt{h}}}
\varphi_h(|[v_h(t)]|)\,d\hs^{N-1}
\le (\bar a+2C')C'.
\end{multline}
Finally for all $h$ and for all $t \in [0,T]$
$$
|D^c v_h(t)|(\Om) \le \frac{C'}{a\sqrt{h}}.
$$
We deduce that \eqref{varbound} holds, and the proof is concluded.
\end{proof}

In order to construct the quasistatic growth of brittle fractures 
in the sense of \cite{FL}
to which $(v_h(t),K_h(t), \gamma_h(t))$ converges, we need the 
following lemma which employees a $\Gamma${-}convergence technique
(see Section \ref{prel}).

\begin{lemma}
\label{gammat}
Let us fix $t \in [0,T]$, and let us consider the functionals
\begin{equation}
\label{basicpblem}
\gs_h(t)(u):=
\int_{\Om} f_h(\nabla u) \,dx+\int_{K_h(t)}\varphi_h(|[u]|)\,d\hs^{N-1}
+a\sqrt{h}|D^cu|(\Om),
\end{equation}
if $u \in BV(\Om)$, $\Sg{g^{\delta_h}(t)}{u} \tsub K_h(t)$, $|[u]| \le
\gamma_h(t)$ on $K_h(t)$,
and $\gs_h(t)(u)=+\infty$ otherwise for $u \in BV(\Om)$.
\par
Let us denote by $\gs(t)$ the $\Gamma${-}limit 
(up to a subsequence) of
$\gs_h(t)$ in the \wstar
topology of $BV(\Om)$. For all $u \in \dom{\gs(t)}$ 
we have $u \in SBV(\Om)$
and $\nabla u \in L^2(\Om;\R^N)$.
Moreover there exists a countable and dense set 
$D \subseteq \dom{\gs(t)}$
such that
setting
\begin{equation}
\label{defgammatlem}
K(t):=\bigcup_{u \in D} \Sg{g(t)}{u}
\end{equation}
we have
\begin{equation}
\label{propgammat}
\Sg{g(t)}{u} \tsub K(t)
\quad\quad \text{ for all }u \in \dom{\gs(t)},
\end{equation}
and
\begin{equation}
\label{lscgammat}
\hs^{N-1}(K(t))  \le \liminf_h \int_{K_h(t)}
\varphi_h(\gamma_h(t))\,d\hs^{N-1}.
\end{equation}
\end{lemma}

\begin{proof}
In order to deal with $\Sg{g}{u}$ as an internal jump,
let us consider $\tilde{\Om} \subseteq \R^N$ open and bounded,
such that $\Omb \subseteq \tilde{\Om}$, and let us set 
$\Om':=\tilde{\Om} \setminus \partial_N \Om$.
Let us consider the following functionals $\gs'_h(t)\,:\,BV(\Om')
\to [0,+\infty]$
\begin{equation}
\label{basicpb}
\gs'_h(t)(u):=
\int_{\Om'} f_h(\nabla u) \,dx+\int_{K_h(t)}\varphi_h(|[u]|)\,d\hs^{N-1}
+a\sqrt{h}|D^cu|(\Om'),
\end{equation}
if $u \in BV(\Om')$, $u=g^{\delta_h}(t)$ on $\Om' \setminus \Om$,
$S(u) \tsub K_h(t)$, $|[u]| \le \gamma_h(t)$ on $K_h(t)$, and
$\gs'_h(t)(u)=+\infty$ otherwise for $u \in BV(\Om')$.
\par
By Proposition \ref{gconvcomp}, up to a subsequence, $\gs'_h(t)$
$\Gamma${-}converges in the \wstar topology of $BV(\Om')$
to a functional $\gs'(t)$. Clearly if $u \in \dom{\gs'(t)}$,
then the restriction of $u$ to $\Om$ belongs
to $\dom{\gs(t)}$. Conversely if $u \in \dom{\gs(t)}$,
the extension of $u$ to $\Om'$ setting $u=g(t)$ on
$\Om' \setminus \Om$ belongs to $\dom{\gs'(t)}$.
Thus we can use $\gs'(t)$ instead of $\gs(t)$.
Let $u \in \dom{\gs'(t)}$: clearly we have $u=g(t)$ on
$\Om' \setminus \Om$.
Moreover since
$$
\gs'_h(t)(u) \ge \int_{\Om'} f_h(\nabla u) \,dx
+\int_{S(u)}\varphi_h(|[u]|)\,d\hs^{N-1}+a\sqrt{h}|D^cu|(\Om'),
$$
if $u \in \dom{\gs'(t)}$ and $u_h \weakst u$ \wlystar in
$BV(\Om')$ with $\gs_h'(u_h) \to \gs'(u)$, by Proposition
\ref{basicgconv} we deduce that $u \in SBV(\Om')$
and
\begin{equation}
\label{compsbv}
\|\nabla u\|^2+\hs^{N-1}(S(u)) \le \gs'(u).
\end{equation}
So we conclude that $\nabla u \in L^2(\Om';\R^N)$, and
$\hs^{N-1}(S(u))<+\infty$.
\par
Let us now consider
$$
\epi{\gs'(t)}:=\{(u,s) \in BV(\Om') \times \R\,:\, \gs'(t)(u) \le s\},
$$
and let $\dst \subseteq \epi{\gs'(t)}$ be countable and dense.
If $\pi\,:\,BV(\Om') \times \R \to BV(\Om')$
denotes the projection on the first factor, let $D:=\pi(\dst)$ and
let us set
\begin{equation}
\label{defgammat}
K(t):= \bigcup_{u \in D} S(u).
\end{equation}
Notice that $K(t)$ is precisely of the form
\eqref{defgammatlem}.
Let us see that $K(t)$
satisfies the properties of the lemma.
\par
Let us prove \eqref{lscgammat}. Let $u_1,\dots,u_k \in D$, and
let $u^h_1,\dots,u^h_k \in BV(\Om')$ be
such that $u^h_i \weakst u_i$ \wlystar in $BV(\Om')$ and
\begin{equation}
\label{gsconv}
\lim_h \gs'_h(t)(u^h_i)=\gs'(t)(u_i),
\quad\quad i=1,\dots,k.
\end{equation}
Setting $u^h:=(u^h_1,\dots,u^h_k)$, by \eqref{gsconv} we have
$$
\sum_{i=1}^k \int_{\Om'} f_h(\nabla u^h_i)\,dx+
\int_{S(u^h)} \varphi_h(|[u^h_1]| \lor \dots \lor
|[u^h_k]|)\,d\hs^{N-1}+a\sqrt{h}|D^cu^h|(\Om')
\le \tilde{C}
$$
with $\tilde{C}$ independent of $h$. By Proposition \ref{basicgconv}
we deduce
\begin{align}
\hs^{N-1} \left( \bigcup_{i=1}^k S(u_i) \right) &\le
\liminf_h \int_{S(u^h)} \varphi_h(|[u^h_1]| \lor \dots \lor |[u^h_k]|)
\,d\hs^{N-1} \\
\nonumber
&\le
\liminf_h \int_{K_h(t)} \varphi_h(\gamma_h(t))\,d\hs^{N-1} \le C',
\end{align}
where $C'$ is given by Corollary \ref{enboundh1lem}.
Taking the sup over all possible $u_1,\dots,u_k$ we get
\begin{equation}
\label{lscgammaht}
\hs^{N-1}(K(t))  \le \liminf_h \int_{K_h(t)}
\varphi_h(\gamma_h(t))\,d\hs^{N-1}
\le C',
\end{equation}
so that \eqref{lscgammat} is proved. In particular we have that
$\hs^{N-1}(K(t))<+\infty$.
\par
Let us come to \eqref{propgammat}. Let $u \in \dom{\gs(t)}$and let us
extend $u$ to $\Om'$ setting
$u=g(t)$ on $\Om' \setminus \Om$. We indicate this extension with $u'$.
We have $u' \in \dom{\gs'(t)}$,
and $S(u')=\Sg{g(t)}{u}$. Let $(u_k,s_k) \in \dst$ be such that $u_k \weakst u'$
\wlystar in
$BV(\Om')$ and $s_k \to \gs'(t)(u')$. By lower semicontinuity of $\gs'(t)$
we have
$$
\gs'(t)(u') \le \liminf_k \gs'(t)(u_k).
$$
Moreover since $\gs'(t)(u_k) \le s_k$, we deduce
$$
\limsup_k \gs'(t)(u_k) \le \lim_k s_k=\gs'(t)(u'),
$$
so that we have $\gs'(t)(u_k) \to \gs'(t)(u')$. By \eqref{compsbv}
we get that $u_k \weak u'$ weakly
in $SBV(\Om')$: since $S(u_k) \tsub K(t)$ for all $k$, and
$\hs^{N-1}(K(t)) <+\infty$, by Ambrosio's
theorem we get $S(u') \tsub K(t)$, i.e. $\Sg{g(t)}{u} \tsub K(t)$.
We conclude that $K(t)$
satisfies \eqref{propgammat}, and the proof is now complete.
\end{proof}

\begin{lemma}
\label{mintlem}
Let $t \in [0,T]$, and let us consider the subsequence of
$(v_h(t),K_h(t),\gamma_h(t))_{h \in \N}$ (which we indicate
with the same symbol) and the rectifiable set $K(t)$
given by Lemma \ref{gammat}.
Then if $v(t) \in SBV(\Om)$ is an accumulation point 
for $(v_h(t))_{h \in \N}$ in
the \wstar topology of $BV(\Om)$,
we have $v(t) \in SBV(\Om)$, $\nabla v(t) \in L^2(\Om;\R^N)$,
\begin{equation}
\label{jumpcontrol}
\Sg{g(t)}{v(t)} \tsub K(t),
\end{equation}
and for all $v \in SBV(\Om)$
\begin{equation}
\label{mint}
\|\nabla v(t)\|^2 \le \|\nabla v\|^2+
\hs^{N-1}(\Sg{g(t)}{v} \setminus K(t)).
\end{equation}
Moreover we have
\begin{equation}
\label{l1conv}
\nabla v_h(t) \weak \nabla v(t)
\quad\quad
\text{weakly in }L^1(\Om;\R^N),
\end{equation}
\begin{equation}
\label{l2conv}
\nabla v_h(t) \1_{E_h(t)} \to \nabla v(t)
\quad\quad
\text{strongly in }L^2(\Om;\R^N),
\end{equation}
where
\begin{equation*}
E_h(t):=
\left\{
x \in \Om\,:\,|\nabla v_h(t)|\le \frac{a\sqrt{h}}{2}
\right\},
\end{equation*}
and
\begin{equation}
\label{ebulk1conv}
\|\nabla v(t)\|^2=\lim_h \int_\Om f_h(\nabla v_h(t))\,dx.
\end{equation}
\end{lemma}

\begin{proof}
Let $v(t)$ be an accumulation point of 
$(v_h(t))_{h \in \N}$ in the
\wstar topology of $BV(\Om)$.
If $\gs(t)$ is the $\Gamma${-}limit of the functional
$\gs_h(t)$ defined in Lemma
\ref{gammat}, by the $\Gamma${-}liminf inequality
and by \eqref{ebound1h} we have
$$
\gs(t)(v(t)) \le \liminf_k \gs_{h_k}(t)(v_{h_k}(t)) \le C',
$$
and so $v(t) \in \dom{\gs(t)}$.
By Lemma \ref{gammat} we conclude $v(t) \in SBV(\Om)$,
$\nabla v(t) \in L^2(\Om;\R^N)$
and $\Sg{g(t)}{v(t)} \tsub K(t)$.
\par
In order to prove \eqref{mint} we follow the Transfer
of Jump of \cite{FL}.
In order to deal with $\Sg{g}{u}$ as an internal jump,
let us consider $\tilde{\Om} \subseteq \R^N$ open and bounded,
and
such that $\Omb \subseteq \tilde{\Om}$. Let us set
$\Om':=\tilde{\Om} \setminus \partial_N \Om$.
By \eqref{defgammatlem} we have
$$
K(t)=\bigcup_{u \in D} S(u),
$$
where $u \in D$ is extended $\Om' \setminus \Om$ setting
$u=g(t)$ on $\Om' \setminus \Om$,
so that $\Sg{g(t)}{u}=S(u)$.
\par
Let $v \in SBV(\Om)$ with $\nabla v \in L^2(\Om;\R^N)$ and
$\hs^{N-1}(\Sg{g(t)}{v})<+\infty$.
Let us consider $w:=v-g(t)$, and let us extend $w$ on $\Om'$
setting $w=0$ on $\Om' \setminus \Om$.
Let $\sigma>0$, and let $u_1,\dots,u_m \in D$ be such that
\begin{equation}
\label{errorsigma}
\hs^{N-1} \left( \Sg{g(t)}{v} \setminus K(t) \right)
< \hs^{N-1} \left( \Sg{g(t)}{v} \setminus \bigcup_{i=1}^m
\Sg{g(t)}{u_i} \right)+\sigma
\end{equation}
Let us fix $G \subseteq \R$ countable and dense:
we recall that for all $r=1,\dots,m$ we have up to a set of
$\hs^{N-1}$-measure zero
$$
S(u_r)= \bigcup_{c_1,c_2 \in G} \partial^* E_{c_1}(u_r) \cap
\partial^*E_{c_2}(u_r),
$$
where $E_c(u_r):=\{x \in \Om'\,:\, x
\text{ is a Lebesgue point for }u_r,\;u_r(x) > c\}$ and
$\partial^*$ denotes the essential boundary
(see \cite[Definition 3.60]{AFP}). Let us orient
$\nu_{u_r}$ in such a way that $u_r^-(x) <u_r^+(x)$ for all
$x \in S(u_r)$, $r=1,\dots,m$, and
let us consider
$$
J_j:=\left\{ x \in \bigcup_{r=1}^m S(u_r) \,:\,
u_i^+(x)-u_i^-(x)>\frac{1}{j} \mbox{ for some } i=1, \ldots, m
\right\},
$$
with $j$ so large that
$$
\hs^{N-1} \left( \bigcup_{r=1}^m S(u_r) \setminus J_j \right)
< \sigma.
$$
Let $U$ be a neighborhood of $\bigcup_{r=1}^m S(u_r)$
such that
\begin{equation}
\label{usmall}
|U|<\frac{\sigma}{j^2},
\quad\quad
\int_U |\nabla w|^2\,dx<\sigma.
\end{equation}
Following \cite[Theorem 2.1]{FL}, we can find a finite
disjoint collection of closed cubes $\{Q_k\}_{k=1, \ldots,n}$
with edge of length $2r_k$, with center $x_k \in S(u_{\tau(k)})$
for some
$\tau(k) \in \{1,\dots,m\}$ and oriented as
the normal $\nu(x_k)$ to $S(u_{\tau(k)})$ at $x_k$, such
that $\bigcup_{k=1}^n Q_k \subseteq U$ and $\hs^{N-1}(J_j
\setminus \bigcup_{k=1}^n Q_k) \le \sigma$.
Let us set
$$
v_k:=u_{\tau(k)}-g(t),
$$
and let $H_k$ denote the intersection of $Q_k$ with the
hyperplane
through $x_k$ orthogonal to $\nu(x_k)$. Following \cite{FL}
we can suppose that
$$
\hs^{N-1} \left( \left[ \bigcup_{r=1}^m S(u_r) \setminus S(v_k)
\right] \cap Q_k \right) < \sigma r_k^{N-1},
$$
and that the following facts hold:
\begin{itemize}
\item[(a)] if $x_k \in \Om$ then $Q_k \subseteq \Om$, and if
$x_k \in \partial_D \Om$, then
$\partial \Om \cap Q_k \subseteq \{y+s\nu(x_k)\,:\,y \in H_k,\,
s \in [-\frac{\sigma r_k}{2},\frac{\sigma r_k}{2}]\}$;
\item[{}]
\item[(b)] $\hs^{N-1} (S(v_k) \cap \partial Q_k)=0$;
\item[{}]
\item[(c)] $r_k^{N-1} < 2 \hs^{N-1} (S(v_k) \cap Q_k)$;
\item[{}]
\item[(d)] $v_k^-(x_k) < c^1_k < c^2_k < v_k^+(x_k)$ and
$c^2_k-c^1_k>\frac{1}{2j}$;
\item[{}]
\item[(e)] $\hs^{N-1} ([S(v_k) \setminus \partial^*
E_{c^s_k}(v_k)] \cap Q_k)<\sigma r_k^{N-1}$ for $s=1,2$;
\item[{}]
\item[(f)] $\hs^{N-1} (\{y \in \partial^*
E_{c^s_k}(v_k) \cap Q_k\,:\,
{\rm dist}(y,H_k) \ge \frac{\sigma}{2} r_k\}) <\sigma r_k^{N-1}$
for $s=1,2$;
\item[{}]
\item[(g)] if $Q^+_k:= \{x\in Q_k \,|\, (x-x_k) \cdot \nu(x_k) >0\}$
and $s=1,2$
\begin{equation}
\label{density1}
\| 1_{E_{c^s_k}(v_k)) \cap Q_k} -1_{Q^+_k}\|_{L^1(\Om')}
<\sigma^2 r_k^N;
\end{equation}
\item[{}]
\item[(h)] $\hs^{N-1} \left( \big( S(w) \setminus S(v_k) \big)
\cap Q_k \right)
< \sigma r_k^{N-1}$ and
$\hs^{N-1} (S(w) \cap \partial Q_k)=0$.
\item[{}]
\end{itemize}
By definition of $D$ in \eqref{defgammatlem}, for all
$k=1,\dots,n$
we can find $v^h_k \in BV(\Om')$ with  $v^h_k=0$ on
$\Om' \setminus \Om$,
$v^h_k \weakst v_k$ \wlystar in $BV(\Om')$,
$S(v^h_k) \tsub K_h(t)$,
$|[v^h_k]| \le \gamma_h(t)$, and such that
\begin{equation}
\label{estvhk}
\limsup_h \int_{\Om'} f_h(\nabla v^h_k) \,dx+
\int_{K_h(t)}\varphi_h(|[v^h_k]|)\,d\hs^{N-1}
+a\sqrt{h}|D^cv^h_k|(\Om') \le \tilde{C},
\end{equation}
with $\tilde{C} \in ]0,+\infty[$. By Proposition
\ref{basicgconv}, we have that $\nabla v^h_k$
converges weakly in $L^1(\Om';\R^N)$, so that we may
assume that $U$ is chosen so that
for $h$ large
\begin{equation}
\label{usmall2}
\sum_{k=1}^n \int_{Q_k} |\nabla v^h_k|\,dx <
\frac{\sigma}{j^2}.
\end{equation}
Let $\eta \in ]0,1[$: we claim that there exists
$\delta>0$ such that for all $k=1,\dots,n$
\begin{equation}
\label{varok1}
\limsup_h |Dv^h_k|(\{0 <|[v^h_k]|< \delta\} \cap Q_k)
\le \eta|Q_k|.
\end{equation}
In fact let $a'<a$ be such that
$$
a's \le \varphi(s)
\quad\quad
\text{for all }s \in [0,1].
$$
Then we have
\begin{align*}
|Dv^h_k| \left( \Big\{ 0 <|[v^h_k]|< \frac{1}{\sqrt{h}} \Big\}
\cap Q_k \right) &=
\int_{\{0 <|[v^h_k]|< \frac{1}{\sqrt{h}}\}\cap Q_k}
|[v^h_k]|\,d\hs^{N-1}
\\
&\le
\frac{1}{a'\sqrt{h}} \int_{\{0 <|[v^h_k]|<
\frac{1}{\sqrt{h}}\} \cap Q_k)} \varphi_h(|[v^h_k]|)\,d\hs^{N-1},
\end{align*}
so that we conclude for $h$ large
\begin{multline}
|Dv^h_k| \left( \Big\{0 <|[v^h_k]|< \delta \Big\} \cap Q_k \right)
\le \frac{1}{a'\sqrt{h}} \int_{\{0 <|[v^h_k]|<
\frac{1}{\sqrt{h}}\} \cap Q_k)} \varphi_h(|[v^h_k]|)\,d\hs^{N-1} \\
+\frac{\delta}{\varphi(1)} 
\int_{\{\frac{1}{\sqrt{h}} <|[v^h_k]|< \delta\} \cap Q_k)}
\varphi_h(|[v^h_k]|)\,d\hs^{N-1}
\le  \left( \frac{1}{a'\sqrt{h}}+\frac{\delta}{\varphi(1)}\right)\tilde{C},
\end{multline}
where $\tilde{C}$ is defined in \eqref{estvhk}.
Taking the limsup in $h$ and choosing $\delta$ small enough,
we have that \eqref{varok1} holds.
\par
Let $\delta$ be as in \eqref{varok1}, and let us set
$$
K^\delta_h(t):=\{x \in K_h(t)\,:\,|[v^h_k]|(x) \ge \delta,\,
\text{ for some }k=1,\dots,n\}.
$$
Then in view of \eqref{usmall2} and of \eqref{varok1},
by the Coarea formula for $BV$ functions 
(see \cite[Theorem 3.40]{AFP}) we have
for $h$ large enough
\begin{multline}
\sum_{k=1}^n \int_{c^1_k}^{c^2_k} \hs^{N-1}
\left( \partial^* E_c(v^h_k) \cap
(Q_k \setminus K^\delta_h(t)) \right) \,dc \le
\sum_{k=1}^n |Dv^h_k|(Q_k \setminus K^\delta_h(t)) \\
=\sum_{k=1}^n \int_{Q_k} |\nabla v^h_k|\,dx+
\sum_{k=1}^n |Dv^h_k|(Q_k \cap \{0 <|[v^h_k]|< \delta\})
\le (1+\eta)\frac{\sigma}{j^2}.
\end{multline}
By the Mean Value Theorem and by property $(d)$ we get that there exist
$c^1_k<c^h_k<c^2_k$ such that
\begin{equation}
\sum_{k=1}^n \hs^{N-1} \left( \partial^* E_{c^h_k}(v^h_k) \cap
(Q_k \setminus K^\delta_h(t)) \right)
\le 2(1+\eta) \frac{\sigma}{j}.
\end{equation}
Following \cite{FL}, by property $(g)$ we have that for $h$ large
$$
\| 1_{E_{c^h_k}(v^h_k) \cap Q_k} -1_{Q^+_k}\|_{L^1(\Om')}
\le \sigma^2 r_k^N.
$$
Then by Fubini's Theorem and by the Mean Value Theorem, we can
find $s^+_k \in [\frac{\sigma r_k}{2},\sigma r_k]$ and
$s^-_k \in [-\sigma r_k, -\frac{\sigma r_k}{2}]$
such that setting $H^+_k:=\{x=y+s^+_k \nu(x_k),\,y \in H_k\}$
and $H^-_k:=\{x=y+s^-_k \nu(x_k),\,y \in H_k\}$
we have
\begin{equation}
\label{estcutk}
\hs^{N-1} \left( H^+_k \setminus
(E_{c^h_k}(v^h_k) \cap Q_k)\right)+
\hs^{N-1} \left( H^-_k \cap (E_{c^h_k}(v^h_k) \cap Q_k) \right)
\le 2\sigma r_k^{N-1}.
\end{equation}
Let $R_k$ be the region between $H^-_k$ and $H^+_k$, i.e.
$$
R_k:=\{x \in Q_k\,:\,x=y+s\nu(x_k),\, y \in H_k,\,s^2_k \le
s \le s^1_k\},
$$
and let us indicate by $R^+_k w$ the reflection in $Q_k$ of
$w_{|Q^+_k \setminus R_k}$ with respect
to $H^+_k$, and by $R^-_k w$ the reflection in $Q_k$ of
$w_{|Q^-_k \setminus R_k}$ with respect
to $H^-_k$.
We can now consider $w_h$ defined in the following way
\begin{equation}
\label{defwk}
w_h:=
\begin{cases}
w &\text{on }\Om' \setminus \bigcup_{k=1}^n R_k \\
R_k^+ w & \text{on }R_k \cap E_{c^h_k}(v^h_k) \\
R_k^- w & \text{on }R_k \setminus E_{c^h_k}(v^h_k).
\end{cases}
\end{equation}
$w_h$ is well defined for $\sigma$ small, and $w_h=0$ on
$\Om' \setminus \Om$. Notice that by construction we have
\begin{equation}
\label{error}
\sum_{k=1}^n \hs^{N-1} \left( \Big( S(w_h) \setminus
K^\delta_h(t) \Big) \cap Q_k \right)
\le e(\sigma),
\end{equation}
where $e(\sigma) \to 0$ as $\sigma \to 0$.
\par
By \eqref{compare1} comparing $v_h(t)$ with
$w_h+g^{\delta_h}(t)$ and in view of
\eqref{errorsigma} and \eqref{error}
we have
\begin{multline}
\label{comparefond}
\int_\Om f_h(\nabla v_h(t))\,dx \le
\int_\Om f_h(\nabla w_h+\nabla g^{\delta_h}(t))\,dx
+\hs^{N-1}(\Sg{g(t)}{w} \setminus K(t)) \\
+\sum_{k=1}^n
\left[
\int_{K^\delta_h(t) \cap Q_k} \varphi_h(w_h \lor
\gamma_h(t))-\varphi_h(\gamma_h(t))
\,d\hs^{N-1}
\right]+e(\sigma).
\end{multline}
Since by construction we have $\gamma_h(t) \ge \delta$ on
$K^\delta_h(t)$, we deduce
$$
\limsup_h \left[
\int_{K^\delta_h(t) \cap Q_k} \varphi_h(w_h \lor
\gamma_h(t))-\varphi_h(\gamma_h(t))
\,d\hs^{N-1}
\right] \le 0.
$$
Moreover we have that
$$
f_h(\nabla w_h+\nabla g^{\delta_h}(t)) \le
|\nabla w_h+\nabla g^{\delta_h}(t)|^2
$$
and by \eqref{usmall} for $h$ large
$$
\|\nabla w_h+\nabla g^{\delta_h}(t)\|^2 \le
\|\nabla w+\nabla g(t)\|^2+e(\sigma).
$$
Since $\sigma$ is arbitrary, we conclude that
\begin{equation}
\label{compfond}
\limsup_h \int_\Om f_h(\nabla v_h(t))\,dx \le
\|\nabla v\|^2+\hs^{N-1}(\Sg{g(t)}{v} \setminus K(t)).
\end{equation}
\par
If $v_{h_m}(t) \weakst v(t)$ \wlystar in $BV(\Om)$,
by Proposition \ref{basicgconv} we have that
$v(t) \in SBV(\Om)$, $\nabla v(t) \in L^2(\Om;\R^N)$,
\begin{equation}
\label{l1convk}
\nabla v_{h_m}(t) \weak \nabla v(t)
\quad\quad
\text{weakly in }L^1(\Om;\R^N),
\end{equation}
and
\begin{equation}
\label{lsc1hm}
\|\nabla v(t)\|^2 \le \liminf_m
\int_{\Om}f_{h_m}(\nabla v_{h_m}(t))\,dx.
\end{equation}
By \eqref{compfond} we obtain
$$
\|\nabla v(t)\|^2 \le \|\nabla v\|^2
+\hs^{N-1}(\Sg{g(t)}{w} \setminus K(t)),
$$
so that \eqref{mint} holds.
\par
Let us now come to \eqref{l1conv}, \eqref{l2conv}
and \eqref{ebulk1conv}.
Since $\Sg{g(t)}{v(t)} \tsub K(t)$ by \eqref{jumpcontrol},
\eqref{mint} implies that
$\nabla v(t)$ is a minimum for the problem
$$
\min \{\|\nabla v\|^2\,:\, v \in SBV(\Om),\, \Sg{g(t)}{v} \tsub K(t)\}.
$$
Since $\nabla v(t)$ is unique by convexity,
by \eqref{l1convk} we deduce that \eqref{l1conv} holds.
Moreover \eqref{ebulk1conv} is a direct consequence of
\eqref{lsc1hm} and \eqref{compfond} with $v=v(t)$.
Finally, notice that $(\nabla v_h(t) \1_{E_h(t)})_{h \in \N}$ 
is bounded in $L^2(\Om;\R^N)$. Since
$\nabla v_h(t) \weak \nabla v(t)$ weakly in $L^1(\Om;\R^N)$ 
and $\nabla v(t) \in L^2(\Om;\R^N)$,
we get $\nabla v_h(t) \1_{E_h(t)} \weak \nabla v(t)$ weakly
in $L^2(\Om;\R^N)$. By \eqref{compfond} with $v=v(t)$
we have
$$
\limsup_h \|\nabla v_h(t) \1_{E_h(t)}\|^2 \le
\limsup_h \int_\Om f_h(\nabla v_h(t))\,dx \le
\|\nabla v(t)\|^2,
$$
so that \eqref{l2conv} holds and the proof is concluded.
\end{proof}

We are now ready to prove Theorem \ref{main1}.

\begin{proof}[Proof of Theorem \ref{main1}]
Point $(a)$ is a consequence of Corollary
\ref{enboundh1lem} and of \eqref{fseub}. Point $(b)$ comes
from Corollary \ref{totvarbound1h}. Let us come to point $(c)$.
\par
Let us consider $B \subseteq [0,T]$ countable and dense.
By a diagonal argument, we may suppose that there exists a unique
subsequence of
$(v_h(t),K_h(t),\gamma_h(t))_{h \in \N}$ 
(which we still denote by the same symbol)
such that Lemma \ref{gammat} holds
for all $t \in B$. For each $t \in B$, let $K(t)$ be the rectifiable
set defined in Lemma \ref{gammat}.
\par
We notice that $K(t)$ is increasing in $t$. In fact if $s<t$
and if $u \in \dom{\gs(s)}$, where $\gs(s)$ is defined in Lemma \ref{gammat},
there exists $u_h \in BV(\Om)$ such that $u_h \weakst u$ \wlystar in
$BV(\Om)$, $\Sg{g^{\delta_h}(s)}{u_h} \tsub K_h(s)$,
$|[u_h]| \le \gamma_h(s)$ and for all $h$
$$
\int_{\Om} f_h(\nabla u_h) \,dx+\int_{K_h(s)}
\varphi_h(|[u_h]|)\,d\hs^{N-1}+c\sqrt{h}|D^cu_h|(\Om)
\le \tilde{C}
$$
for some $\tilde{C}$ independent of $h$.
Let us set $v_h:=u_h-g^{\delta_h}(s)+g^{\delta_h}(t)$.
Since $K_h(s) \tsub K_h(t)$ and $\gamma_h(s) \le \gamma_h(t)$,
we have that
$\Sg{g^{\delta_h}(t)}{v_h} \tsub K_h(t)$, $|[v_h]| \le \gamma_h(t)$;
moreover for all $h$
$$
\int_{\Om} f_h(\nabla v_h) \,dx+\int_{K_h(t)}\varphi_h(|[v_h]|)\,d\hs^{N-1}
+a\sqrt{h}|D^cv_h|(\Om)
\le \tilde{C}'
$$
with $\tilde{C}'$ independent of $h$.
We deduce that $u-g(s)+g(t) \in \dom{\gs(t)}$.
Then by definition \eqref{defgammatlem} and by
\eqref{propgammat}, we obtain $K(s) \tsub K(t)$.
\par
Since $\{t \to K(t): t \in [0,T]\}$ is increasing, setting
$$
K^-(t):=\bigcup_{s \in B, s \le t} K(s),
\quad\quad
K^+(t):=\bigcap_{s \in B, s \ge t} K(s),
$$
there exists a countable set $B' \subseteq [0,T] \setminus B$
such that we have $K^-(t)=K^+(t)$
for all $t \in [0,T] \setminus B'$. For all such $t$'s
let us set $K(t):=K^-(t)=K^+(t)$.
\par
Clearly Lemma \ref{gammat} and Lemma \ref{mintlem} hold for all
$t \in [0,T] \setminus (B \cup B')$. In fact, up to a further
subsequence, we may apply
Lemma \ref{gammat} obtaining $\tilde{K}(t)$ with the required
properties and such
that $K(s_1) \tsub \tilde{K}(t) \tsub K(s_2)$ for all $s_1,s_2 \in B$
and $s_1<t<s_2$.
Then we get $\tilde{K}(t)=K(t)$.
\par
Up to a further subsequence relative to the elements of $B'$,
we find
$$
\{t \to K(t),\,t \in [0,T] \}
$$
such that Lemma \ref{gammat} and Lemma \ref{mintlem} hold for every
$t \in [0,T]$. Notice that
in particular $\hs^{N-1}(K(t)) \le C'$, where $C'$ is given by
\eqref{ebound1h}.
\par
Let $v(t)$ be a minimum for the following problem
\begin{equation}
\label{defvt}
\min \left\{\ \|\nabla v\|^2\,:\,v \in SBV(\Om),\,
\Sg{g(t)}{v} \tsub K(t) \right\}.
\end{equation}
Notice that problem \eqref{defvt} is well posed since $K(t)$ has
finite $\hs^{N-1}${-}measure: moreover
by strict convexity we have that $\nabla v(t)$ is uniquely determined.
\par
We claim that $\{t \to (v(t),K(t)),\,t \in [0,T] \}$
is a quasistatic growth of brittle fractures
in the sense of \cite{FL}, that is in the sense of Theorem \ref{qse}.
In fact $K(t)$ is increasing and satisfies the unilateral minimality
property \eqref{mint} by construction.
In order to prove the claim, we have just to prove the
{\it nondissipativity} condition
\begin{equation}
\label{nondiss}
\Es(t)=\Es(0)+2\int_0^t (\nabla v(\tau),\nabla
\dot{g}(\tau))_{L^2(\Om;\R^N)}\,d\tau,
\end{equation}
where $\Es(t):=\|\nabla v(t)\|^2+\hs^{N-1}(K(t))$
for all $t \in [0,T]$.
First of all for all $t \in [0,T]$ we have
\begin{equation}
\label{energybelow}
\Es(t) \ge \Es(0)+2 \int_0^t (\nabla u(\tau),\nabla
\dot{g}(\tau))_{L^2(\Om;\R^N)}\,d\tau.
\end{equation}
In fact as noticed in \cite{GP}, using the minimality
property \eqref{mint},
the map $\{t \to \nabla v(t)\}$ is continuous at all the
continuity points of
$\{t \to \hs^{N-1}(K(t))\}$, in particular it is continuous
up to a countable set in $[0,T]$. Given $t \in [0,T]$ and $k>0$,
let us set
$$
s_i^k:= \frac{i}{k}t,
\quad\quad
v^k(s):=v(s_{i+1}^k) \quad \text{ for }s_i^k< s \le s_{i+1}^k
$$
for all $i=0,1,\dots,k$. By \eqref{mint}, comparing
$v(s^k_i)$ with $v(s^k_{i+1})-g(s^k_{i+1})+g(s^k_{i})$,
it is easy to see that
\begin{equation*}
\Es(t) \ge \Es(0)+
2\int_0^t (\nabla v^k(\tau),
\nabla \dot{g}(\tau))_{L^2(\Om;\R^N)} \,d\tau +e(k),
\end{equation*}
where $e(k) \to 0$ as $k \to +\infty$.
By the continuity property of $\nabla v$,
passing to the limit for $k \to +\infty$ we deduce that
\eqref{energybelow} holds.
On the other hand for all $t \in [0,T]$ we have that
\begin{equation}
\label{energyabove1}
\Es(t) \le \Es(0)+2 \int_0^t (\nabla u(\tau),
\nabla \dot{g}(\tau))_{L^2(\Om;\R^N)}\,d\tau.
\end{equation}
In fact by Proposition \ref{gammaconvzero} and by property
\eqref{comparezero1} we get
\begin{equation}
\label{conv0}
\fs_h(0,v_h(0)) \to \Es(0).
\end{equation}
Moreover by Lemma \ref{mintlem} we have that for all $t \in [0,T]$
$$
\nabla v_h(t) \1_{E_h(t)} \to \nabla v(t)
\quad\quad
\text{strongly in }L^2(\Om;\R^N),
$$
where
$$
E_h(t):=
\left\{
x \in \Om\,:\, |\nabla v_h(t)| \le \frac{a\sqrt{h}}{2}
\right\}.
$$
By \eqref{en1above} and by the very definition of $f_h$
we deduce
\begin{align}
\label{en1abovebis}
\fs_h(t,v_h(t)) \le &\fs_h(0,v_h(0))
+2 \int_0^t (\nabla v_h(\tau) \1_{E_h(\tau)},
\nabla \dot{g}(\tau))_{L^2(\Om;\R^N)}\,d\tau \\
\nonumber
&+a\sqrt{h} \int_0^t \int_{\Om \setminus E_h(\tau)}
|\nabla \dot{g}(\tau)|\,dx\,d\tau+e(h)
\end{align}
where $e(h) \to 0$ as $h \to +\infty$.
Notice that by \eqref{ebound1h} we have
$$
\frac{a}{2}h|\Om \setminus E_h(\tau)| \le \sqrt{h}
\int_{\Om \setminus E_h(\tau)} |\nabla v_h(\tau)|\,dx
\le \frac{2}{a}
\int_{\Om \setminus E_h(\tau)} f_h(\nabla v_h(\tau))\,dx
\le \frac{2}{a}C'.
$$
We deduce that
\begin{align}
\label{nullerror}
\sqrt{h}\int_{\Om \setminus E_h(\tau)} |\nabla \dot{g}(\tau)|\,dx
&\le
\left( \int_{\Om \setminus E_h(\tau)}
|\nabla \dot{g}(\tau)|^2\,dx \right)^{\frac{1}{2}}
\sqrt{h|\Om \setminus E_h(\tau)|} \\
\nonumber
& \le \frac{2\sqrt{C'}}{a}\left( \int_{\Om \setminus E_h(\tau)}
|\nabla \dot{g}(\tau)|^2\,dx \right)^{\frac{1}{2}} \to 0
\end{align}
uniformly in $\tau$ as $h \to +\infty$ by equicontinuity of
$\nabla \dot{g}(\tau)$. Then passing to the
limit for $h \to +\infty$ in \eqref{en1abovebis},
in view of \eqref{ebulk1conv}, \eqref{lscgammat}, \eqref{l2conv},
\eqref{conv0} and \eqref{nullerror}
we deduce that \eqref{energyabove1}
holds. This proves that
\eqref{nondiss} holds, and so $\{t \to (v(t),K(t))\,:\,t \in [0,T]\}$ 
is a quasistatic
growth of brittle fractures in the sense of \cite{FL}.
\par
In order to conclude the proof, let us see that
\eqref{gradconv1}, \eqref{totconv1},
\eqref{bulkconv1} and \eqref{fracconv1} hold.
By \eqref{en1abovebis} we deduce that for all $t \in [0,T]$
\begin{equation}
\label{totconv2bis}
\fs_h(t,v_h(t)) \to \Es(t),
\end{equation}
so that by \eqref{ebulk1conv} and \eqref{lscgammat}
we deduce that
\begin{equation}
\label{jumpconv2}
\hs^{N-1}(K(t))=
\lim_h \int_{K_h(t)} \varphi_h(\gamma_h(t))\,d\hs^{N-1},
\quad\quad
a\sqrt{h}|D^cv_h(t)|(\Om) \to 0.
\end{equation}
Theorem \ref{main1} is now completely proved in view of
the rescaling
\eqref{riscaledisp1}, of \eqref{defgammapiccoloh},
\eqref{defvarphi1h} and \eqref{fseub}.
\end{proof}


\section{Proof of Theorem \ref{main2}}
\label{main2sec}
In this section we will give the proof of Theorem \ref{main2}.
Let $\{t \to (u_h(t),\Gamma_h(t),\psi_h(t))\,:\,t \in [0,T]\}$
be the piecewise constant interpolation
of a discrete in time evolution of cohesive fracture in $\Om_h$ relative
the subdivision
$I_{\delta_h}:=\{0=t^{\delta_h}_0<\dots<t^{\delta_h}_{N_{\delta_h}}=T\}$,
the preexisting crack configuration $({\bar \Gamma}_h,{\bar \psi}_h)$ given by
\eqref{initconfig}
and the boundary displacement $h^\alpha g(t,\frac{x}{h})$ with
$\alpha \in ]0,\frac{1}{2}[$.
\par
In order to prove Theorem \ref{main2}, proceeding as in Section
\ref{main1sec},
it is convenient to rescale $u_h$ and $\Gamma_h$ in the
following way: for all
$t \in [0,T]$ let $v_h(t) \in BV(\Om)$ and
$K_h(t) \tsub \Om \cup \partial_D \Om$ be given by
\begin{equation}
\label{riscaledisp2}
v_h(t,x):=\frac{1}{h^\alpha}u_h \left( t,hx \right),
\quad\quad
K_h(t):=\frac{1}{h} \Gamma_h(t).
\end{equation}
Let us moreover set
\begin{equation}
\label{defgammapiccoloh2}
\gamma_h(t,x):=\frac{1}{h^\alpha} \psi_h(t,hx)=
\max_{s \le t} |[v_h(s)](t,x)| \lor \bar \gamma(x)
\quad\quad
t \in [0,T], x \in \Om.
\end{equation}
It turns that $\{t \to (v_h(t),K_h(t),\gamma_h(t))\,:\,t \in [0,T]\}$
is the piecewise constant interpolation
of a discrete in time evolution of cohesive fractures in $\Om$
relative to the subdivision $I_{\delta_h}$,
the preexisting crack configuration $(\bar \Gamma,\bar \gamma)$
and boundary displacement $g(t)$ with
respect to the basic total energy
\begin{equation}
\label{riscenergy2}
\int_\Om f_h(\nabla v) \,dx+
h^{1-2\alpha} \int_{\Sg{g^{\delta_h}(t)}{v} \cup K_h(t)}
\varphi_{h}(|[v]| \lor \gamma_h(t))\,d\hs^{N-1}
+ah^{1-\alpha}|D^cv|(\Om),
\end{equation}
where
\begin{equation}
\label{defvarphi2h}
\varphi_h(s):=\varphi(h^\alpha s),
\end{equation}
and
\begin{equation}
\label{deff2h}
f_h(\xi):=
\begin{cases}
|\xi|^2 & \text{if } |\xi| \le \frac{a h^{1-\alpha}}{2}
\\ \\
\frac{a^2h^{2(1-\alpha)}}{4}+a h^{1-\alpha}(|\xi|-\frac{ah^{1-\alpha}}{2})
&\text{if } |\xi| \ge \frac{a h^{1-\alpha}}{2}.
\end{cases}
\end{equation}
We have that the following facts hold:
\begin{itemize}
\item[]
\item[(a)] for all $t \in [0,T]$
\begin{equation}
\label{infty2}
\|v_h(t)\|_\infty \le \|g^{\delta_h}(t)\|_\infty \le C;
\end{equation}
\item[]
\item[(b)] for all $w \in BV(\Om)$ we have
\begin{multline}
\label{comparezero2}
\int_{\Om} f_h(\nabla v_h(0)) \,dx
+h^{1-2\alpha} \int_{K_h(0)} \varphi_h(|[v_h(0)]| \lor \bar \gamma)
\,d\hs^{N-1} +ah^{1-\alpha}|D^c v_h(t)|(\Om) \\
\le \int_{\Om} f_h(\nabla w) \,dx
+h^{1-2\alpha} \int_{\Sg{g^{\delta_h}(0)}{w} \cup \bar \Gamma}
\varphi_h(|[w]| \lor \bar \gamma)\,d\hs^{N-1} 
+a h^{1-\alpha}|D^cw|(\Om);
\end{multline}
\item[]
\item[(c)] for all $w \in BV(\Om)$ we have
\begin{multline}
\label{compare2}
\int_{\Om} f_h(\nabla v_h(t)) \,dx
+h^{1-2\alpha} \int_{K_h(t)} \varphi_h(\gamma_h(t))\,d\hs^{N-1} 
+ah^{1-\alpha}|D^c v_h(t)|(\Om) \\
\le \int_{\Om} f_h(\nabla w) \,dx
+h^{1-2\alpha} \int_{\Sg{g^{\delta_h}(t)}{w} \cup K_h(t)}
\varphi_h(|[w]| \lor \gamma_h(t)\,d\hs^{N-1} 
+a h^{1-\alpha}|D^cw|(\Om).
\end{multline}
\end{itemize}
Let us set for all $v \in BV(\Om)$ and for all $t \in [0,T]$
\begin{equation}
\label{deffs2}
\fs_h(t,w):=\int_{\Om} f_h(\nabla w) \,dx
+h^{1-2\alpha}\int_{\Sg{g^{\delta_h}(t)}{w} \cup K_h(t)}
\varphi_h(|[w]| \lor \gamma_h(t))\,d\hs^{N-1} 
+a h^{1-\alpha}|D^cw|(\Om).
\end{equation}
Notice that
\begin{equation}
\label{fseub2}
\fs_h(t,v_h(t))=\frac{1}{h^{N+2\alpha-2}} \Es(t,u_h(t)),
\end{equation}
where $\Es(t,u_h(t))$ is defined in \eqref{energydefth}.
We can now prove Theorem \ref{main2}.

\begin{proof}[Proof of Theorem \ref{main2}] 
By Lemma \ref{energyabove}
we obtain for all $t \in [0,T]$
\begin{equation}
\label{esttoten2}
\fs_h(t,v_h(t)) \le \fs_h(0,v_h(0))+
\int_{0}^{t_h} \int_\Om 
f'_h(\nabla v_h(\tau)) \nabla \dot{g}(\tau)\,dx\,d\tau+e(h),
\end{equation}
where $e(h) \to 0$ as $h \to +\infty$, and
$t_h:=t^{\delta_h}_{i_h}$ is the step
discretization point of $I_{\delta_h}$ such that
$t^{\delta_h}_{i_h} \le t<t^{\delta_h}_{i_h+1}$.
By \eqref{comparezero2} comparing $v_h(0)$ and $g(0)$ we have
\begin{multline}
\label{energybound2tris}
\int_\Om f_h(\nabla v_h(0)) \,dx+
h^{1-2\alpha}
\left[
\int_{\Sg{g(0)}{v_h(0)} \cup \bar \Gamma} \varphi_{h}(|[v_h(0)]|
\lor \bar \gamma)-\varphi_{h}(\bar \gamma)\,d\hs^{N-1}
\right]
\\+a h^{1-\alpha}|D^cv_h(0)|(\Om) \le \|\nabla g(0)\|^2.
\end{multline}
By \eqref{compare2} comparing $v_h(t)$ and $g^{\delta_h}(t)$
we obtain
\begin{equation}
\label{energybound2bis}
\int_\Om f_h(\nabla v_h(t))\,dx \le
\|\nabla g^{\delta_h}(t)\|^2,
\end{equation}
and since we have
$$
\int_\Om |f'_h(\nabla v_h(\tau))|^2\,dx \le
4\int_\Om f_h(\nabla v_h(\tau))\,dx ,
$$
by \eqref{esttoten2} we deduce that
\begin{equation}
\label{ebound2}
\int_\Om f_h(\nabla v_h(t)) \,dx+
h^{1-2\alpha} \int_{K_h(t) \setminus \bar \Gamma}
\varphi_{h}(\gamma_h(t))\,d\hs^{N-1}
+a h^{1-\alpha}|D^cv|(\Om) \le C'
\end{equation}
with $C'$ independent of $h$ and of $t$.
By \eqref{gbarfinite} and \eqref{infty2} and
following Corollary \ref{totvarbound1h}
we deduce that $(v_h(t))_{h \in \N}$ 
is bounded in $BV(\Om)$, and this proves point $(a)$.
\par
Let $v(t)$ be an accumulation point for 
$(v_h(t))_{h \in \N}$
in the \wstar topology of $BV(\Om)$,
and let us consider $\tilde{\Om} \subseteq \R^N$
open and bounded, and
such that $\Omb \subseteq \tilde{\Om}$.
Let us set $\Om':=\tilde{\Om} \setminus \partial_N \Om$.
Then we can extend $v_h(t)$ and $v(t)$ to $\Om'$
setting $v_h(t)=g^{\delta_h}(t)$ and
$v(t)=g(t)$ on $\Om' \setminus \Om$ respectively. We have
$v_{h_j}(t) \weakst v(t)$ \wlystar in $BV(\Om')$
for a suitable $h_j \nearrow +\infty$, and
\begin{equation}
\label{coerc2}
\int_{\Om'} f_{h_j}(\nabla v_{h_j}(t)) \,dx+
h_j^{1-2\alpha} \int_{S(v_{h_j}(t)) \setminus \bar \Gamma}
\varphi_{h_j}(|[v_{h_j}(t)]|)\,d\hs^{N-1}
+a h_j^{1-\alpha}|D^cv_{h_j}(t)|(\Om') \le \tilde{C}
\end{equation}
with $\tilde{C}$ independent of $j$. In particular we have
$$
\int_{\Om'} f_{h_j}(\nabla v_{h_j}(t)) \,dx+
\int_{S(v_{h_j}(t))}
\varphi_{h_j}(|[v_{h_j}(t)]|)\,d\hs^{N-1}
+a h_j^{1-\alpha}|D^cv_{h_j}(t)|(\Om') \le \tilde{C}'
$$
with $\tilde{C}'$ independent of $j$.
Then by Proposition \ref{basicgconv} we have that 
$v(t) \in SBV(\Om)$,
\begin{equation}
\label{weak2}
\nabla v_{h_j}(t) \weak \nabla v(t)
\quad\quad
\text{weakly in }L^1(\Om;\R^N),
\end{equation}
and
\begin{equation}
\label{weaklsc2}
\|\nabla v(t)\|^2 \le 
\liminf_j \int_{\Om} f_{h_j}(\nabla v_{h_j}(t)) \,dx.
\end{equation}
Finally, if we consider for all Borel sets $B \subseteq \Om'$
$$
\lambda_j(B):=\int_{B \cap S(v_{h_j}(t))}
\varphi_{h_j}(|[v_{h_j}(t)]|)\,d\hs^{N-1}
$$
and if (up to a subsequence) $\lambda_j \weakst \lambda$
\wlystar in the sense of measures,
we deduce following Proposition \ref{basicgconv} that
$$
\hs^{N-1} \res S(v(t)) \le \lambda
\quad\quad
\text{as measures.}
$$
Since by \eqref{coerc2} we have 
$\lambda(\Om' \setminus \bar \Gamma)=0$, then we have
$S(v(t)) \tsub \bar \Gamma$, that is 
$\Sg{g(t)}{v(t)} \tsub \bar \Gamma$.
Point $(b)$ is thus proved.
\par
Let us come to point $(c)$. Let us suppose that 
$\varphi(s)=1$ for $s \ge \bar s$,
and $\bar \gamma \ge \eps>0$. Let us consider $v \in SBV(\Om)$ with
$\Sg{g(t)}{v} \tsub \bar \Gamma$. Comparing $v_h(t)$ 
with $v-g(t)+g^{\delta_h}(t)$
by minimality property \eqref{compare2} we obtain
\begin{multline}
\int_{\Om} f_h(\nabla v_h(t)) \,dx
+h^{1-2\alpha} \int_{K_h(t)} \varphi_h(\gamma_h(t))\,d\hs^{N-1} 
+ah^{1-\alpha}|D^c v_h(t)|(\Om) \\
\le \int_{\Om} f_h(\nabla v-\nabla g(t)+\nabla g^{\delta_h}(t)) \,dx
+h^{1-2\alpha} \int_{K_h(t)}
\varphi_h(|[v]| \lor \gamma_h(t))\,d\hs^{N-1}.
\end{multline}
Since for $h$ large we have
$$
\varphi_h(\gamma_h(t)) =\varphi_h(|[v]| \lor \gamma_h(t)),
$$
we deduce that
\begin{align}
\label{compare2h}
\int_{\Om} f_h(\nabla v_h(t)) \,dx
+ah^{1-\alpha}|D^c v_h(t)|(\Om)
&\le \int_{\Om} f_h(\nabla v-\nabla g(t)+\nabla g^{\delta_h}(t)) \,dx
\\ \nonumber &\le
\|\nabla v-\nabla g(t)+\nabla g^{\delta_h}(t)\|^2.
\end{align}
Let $h_j \nearrow +\infty$ such that 
$v_{h_j}(t) \weakst v(t)$ \wlystar in $BV(\Om)$.
By \eqref{weak2} we have that 
$\nabla v_{h_j}(t) \weak \nabla v(t)$
weakly in $L^1(\Om;\R^N)$ and in view of \eqref{weaklsc2} 
we deduce that
$$
\|\nabla v(t)\|^2 \le \|\nabla v\|^2,
$$
so that $v(t)$ is a minimizer of
$$
\min \{\|\nabla v\|^2\,:\, v \in SBV(\Om),\,
\Sg{g(t)}{v} \tsub \bar \Gamma\}.
$$
By strict convexity, we have that $\nabla v(t)$ 
is uniquely determined, and so we deduce that
$\nabla u_h(t) \weak \nabla u(t)$ weakly in $L^1(\Om;\R^N)$.
Choosing $v=u(t)$ in \eqref{compare2h} and taking the limsup in 
$h$ we have
$$
\limsup_h \int_{\Om} f_h(\nabla v_h(t)) \,dx \le \|\nabla u(t)\|^2,
$$
so that
$$
\lim_h \int_{\Om} f_h(\nabla v_h(t)) \,dx = \|\nabla u(t)\|^2.
$$
The proof is concluded thank to the rescaling \eqref{riscaledisp2}.
\end{proof}

\section{Proof of Theorem \ref{main3}}
\label{main3sec}
In this section we will give the proof of Theorem \ref{main3}.
Let $\{t \to (u_h(t),\Gamma_h(t),\psi_h(t))\,:\,t \in [0,T]\}$
be the piecewise constant interpolation
of a discrete in time evolution of cohesive fracture in $\Om_h$ 
relative the subdivision
$I_{\delta_h}:=\{0=t^{\delta_h}_0<\dots<t^{\delta_h}_{N_{\delta_h}}=T\}$, 
the preexisting crack configuration
$({\bar \Gamma}_h,{\bar \psi}_h)$ given by \eqref{initconfig}
and the boundary displacement $h^\alpha g(t,\frac{x}{h})$ with
$\alpha \ge \frac{1}{2}$.
\par
In order to prove Theorem \ref{main3}, proceeding as in Section
\ref{main1sec},
it is convenient to rescale $u_h$ and $\Gamma_h$ in the following way:
for all
$t \in [0,T]$ let $v_h(t) \in BV(\Om)$ and
$K_h(t) \tsub \Om \cup \partial_D \Om$
be given by
\begin{equation}
\label{riscaledisp3}
v_h(t,x):=\frac{1}{h^\alpha}u_h \left( t,hx \right),
\quad\quad
K_h(t):=\frac{1}{h} \Gamma_h(t).
\end{equation}
Let us moreover set
\begin{equation}
\label{defgammapiccoloh3}
\gamma_h(t,x):=\frac{1}{h^\alpha} \psi_h(t,hx)=
\max_{s \le t} |[v_h(s)](t,x)| \lor \bar \gamma(x),
\quad\quad
t \in [0,T], x \in \Om.
\end{equation}
It turns out that $\{t \to (v_h(t),K_h(t),\gamma_h(t))\,:\, t\ in [0,T]\}$ 
is the piecewise constant interpolation of 
a discrete in time evolution of cohesive fractures in $\Om$
relative to the preexisting crack configuration
$(\bar \Gamma,\bar \gamma)$
and boundary displacement $g(t)$ with
respect to the basic total energy
\begin{equation}
\label{riscenergy3}
h^{2\alpha-1}\int_\Om f_h(\nabla v) \,dx+
\int_{\Sg{g^{\delta_h}}{v}}
\varphi_{h}(|[v]| \lor \gamma_h(t))\,d\hs^{N-1}+a h^\alpha|D^cv|(\Om),
\end{equation}
where
\begin{equation}
\label{defvarphi3h}
\varphi_h(s):=\varphi(h^\alpha s),
\end{equation}
and
\begin{equation}
\label{deff3h}
f_h(\xi):=
\begin{cases}
|\xi|^2 & \text{if } |\xi| \le \frac{a h^{1-\alpha}}{2}
\\ \\
\frac{a^2h^{2(1-\alpha)}}{4}
+a h^{1-\alpha}(|\xi|-\frac{ah^{1-\alpha}}{2})
&\text{if } |\xi| \ge \frac{a h^{1-\alpha}}{2}.
\end{cases}
\end{equation}
Notice that by Proposition \ref{discrevol} we have
\begin{equation}
\label{infty3}
\|v_h(0)\|_\infty \le \|g(0)\|_\infty \le C,
\end{equation}
and for all $w \in BV(\Om)$ we have
\begin{multline}
\label{comparezero3}
h^{2\alpha-1} \int_{\Om} f_h(\nabla v_h(0)) \,dx
+\int_{\Sg{g^{\delta_h}(0)}{v_h(0)} \cup \bar \Gamma}
\varphi_h(|[v_h(0)]| \lor \bar \gamma)\,d\hs^{N-1}
+ah^\alpha|D^c v_h(0)|(\Om) \\
\le h^{2\alpha-1} \int_{\Om} f_h(\nabla w) \,dx
+\int_{\Sg{g^{\delta_h}(0)}{w} \cup \bar \Gamma}
\varphi_h(|[w]| \lor \bar \gamma)\,d\hs^{N-1} +a h^\alpha|D^cw|(\Om).
\end{multline}
We can now prove Theorem \ref{main3}.

\begin{proof}[Proof of Theorem \ref{main3}]
Comparing $v_h(0)$ and $w=-C$ by means of \eqref{comparezero3} we have
\begin{multline}
\label{energybound3tris}
h^{2\alpha-1} \int_{\Om} f_h(\nabla v_h(0)) \,dx
+\int_{\Sg{g^{\delta_h}(0)}{v_h(0)} \cup \bar \Gamma}
\varphi_h(|[v_h(0)]| \lor \bar \gamma)\,d\hs^{N-1}
+ah^\alpha|D^c v_h(0)|(\Om) \\
\le \hs^{N-1}(\bar{\Gamma} \cup \partial_D \Om).
\end{multline}
As a consequence, following Corollary \ref{totvarbound1h}, we obtain
$$
|Dv_h(0)|(\Om) \le C'
$$
with $C'$ independent of $h$. Since moreover $\|v_h(0)\|_\infty \le C$ by
\eqref{infty3}, we deduce that $(v_h(0))_{h \in \N}$ is bounded in $BV(\Om)$.
Let $v$ be an
accumulation point for $(v_h(0))_{h \in \N}$ 
in the \wstar topology of $BV(\Om)$.
Let us prove that
$v \in SBV(\Om)$ and that $\nabla v=0$: in fact we have that for all
$\xi \in \R^N$
$$
\tilde{f}_h(\xi) \le h^{2\alpha-1}f_h(\xi)
$$
where
\begin{equation}
\label{defftilde3h}
\tilde{f}_h(\xi):=
\begin{cases}
|\xi|^2 & \text{if } |\xi| \le \frac{a h^\alpha}{2}
\\ \\
\frac{a^2h^{2\alpha}}{4}+a h^\alpha(|\xi|-\frac{ah^\alpha}{2})
&\text{if } |\xi| \ge \frac{a h^\alpha}{2}.
\end{cases}
\end{equation}
We deduce that there exists $C''$ independent of $h$ such that
for all $h$
\begin{equation}
\label{bound3}
\int_{\Om} \tilde{f}_h(\nabla v_h(0))\,dx+
\int_{S(v_h(0))} \varphi_h(|[v_h(0)]|)\,d\hs^{N-1}
+ch^\alpha|D^cv_h(0)|(\Om) \le C''.
\end{equation}
By Proposition \ref{basicgconv}, we obtain that $v \in SBV(\Om)$
and that
$\nabla v_h(0) \weak \nabla v$ weakly in $L^1(\Om;\R^N)$.
By \eqref{energybound3tris} we obtain that
$$
\|\nabla v_h(0)\|_{L^1(\Om;\R^N)} \le
\frac{\hs^{N-1}(\bar \Gamma \cup \partial_D \Om)+1}{ah^\alpha},
$$
so that we deduce $\nabla v=0$, that is $v$ is piecewise
constant in $\Om$.
Finally taking the limit in \eqref{comparezero3} with 
$w$ piecewise constant, then we get exactly
\eqref{minpiece}, so that the proof of Theorem \ref{main3} 
is concluded.
\end{proof}

\section{A relaxation result}
\label{relres}
In this section, we prove a relaxation result we used
in order to study the discrete in time evolution of fractures 
in the cohesive case.
\par
Let $f:\R \to [0,+\infty[$ be convex, $f(0)=0$ and with superlinear
growth, i.e.
$$
\limsup_{|\xi| \to +\infty} \frac{f(\xi)}{|\xi|}=+\infty.
$$
Let $\varphi:[0,+\infty[ \to [0,+\infty[$ be increasing, concave,
and such that $\varphi(0)=0$.
Notice that if $a:=\varphi'(0)<+\infty$, we have
\begin{equation}
\label{phiabove}
\varphi(s) \le as
\quad \text{for all }s \in[0,+\infty[.
\end{equation}
Let $\Om$ be a Lipschitz bounded open set in $\R^N$, and let
$\partial_D \Om \subseteq \partial \Om$
be open in the relative topology.
Let $\Gamma$ be a rectifiable set in $\Om \cup \partial_D \Om$,
and let $\psi$ be a positive function
defined on $\Gamma$. Let us extend $\psi$ to
$\Om \cup \partial_D \Om$ setting $\psi=0$ outside $\Gamma$.
Let $g \in W^{1,1}(\Om)$: we may assume that $g$ is
extended to the whole $\R^N$, and
we indicate this extension still by $g$.
\par
We will study the following functional
\begin{equation}
\label{func}
F(u):=
\begin{cases}
\int_\Om f(|\nabla u|) \,dx+ \int_{\Sg{g}{u} \cup \Gamma}
\varphi(|[u]| \lor \psi)\,d\hs^{N-1} \\ \\
+\infty
\end{cases}
\begin{array}{l}
\text{if }u \in SBV(\Om) \\ \\ \text{otherwise in }BV(\Om),
\end{array}
\end{equation}
where $\Sg{g}{u}$ is defined in \eqref{sgjump}, and
$a \lor b:=\max\{a,b\}$ for all $a,b \in \R$.
The functional \eqref{func} naturally appears
(see Section \ref{discrevolsec})
when dealing with quasistatic growth of fractures in the cohesive case,
 where one is required to look for its minima.
We are led to compute the relaxation of $F$ with respect to the
strong topology of $L^1(\Om)$.
The relaxation in the case $\Gamma=\emptyset$ (without boundary
conditions but
without superlinear growth on $f$) has been proved in \cite{BBB}.
Let
\begin{equation}
\label{bulkrelax}
f_1(\xi):= \inf\{f(\xi_1)+a|\xi_2|\,:\, \xi_1+\xi_2=\xi\},
\end{equation}
where $a:=\varphi'(0)$. We have that the following result holds.

\begin{proposition}
\label{relaxation}
The relaxation of the functional \eqref{func} with respect to the
\wstar topology of $BV(\Om)$
is given by $\overline{F}\,:\,BV(\Om) \to [0,+\infty]$ defined as
\begin{equation}
\label{fbar}
\overline{F}(u):=
\int_\Om f_1(|\nabla u|) \,dx+\int_{\Sg{g}{u} \cup \Gamma}
\varphi(|[u]|\lor \psi)\,d\hs^{N-1}
+a|D^cu|,
\end{equation}
where $a=\varphi'(0)$ and $f_1$ is defined in \eqref{bulkrelax}.
\end{proposition}

In order to prove Proposition \ref{relaxation}, we need some
preliminaries.
\par\noindent
Let $I \subseteq \R$ be a finite union of disjoint intervals, and
let $J\subseteq I$ be a countable set.
Let us consider the functional
\begin{equation}
\label{funcmeas}
\fs(\mu):=\int_I f_1(|\phi_{\mu}|)\,dx+
\sum_{t \in S_\mu \setminus J} \varphi(|\mu(\{t\})|)
+\sum_{t \in J} \varphi(|\mu(\{t\})| \lor \psi(t))+a|\mu^c|(I)
\end{equation}
defined for all $\mu \in \ms_b(I;\R^k)$, i.e. $\mu$ is a bounded
$\R^k$-valued Radon measure on $I$.
Here $\phi_{\mu}$ is the density of the absolutely continuous part
$\mu^a$ of $\mu$,
$S_\mu$ is the set of atoms of $\mu$,
$\mu^c:=\mu-\mu^a-\mu \res S_\mu$, $\psi$ is a strictly positive
function defined on $J$, $a=\varphi'(0)$ and
$f_1$ is defined in \eqref{bulkrelax}.

\begin{lemma}
\label{lscmeasure}
The functional $\fs$ defined in \eqref{funcmeas} is lower semicontinuous
with respect to the $\text{weak}^*$
convergence in the sense of measures.
\end{lemma}

\begin{proof}
Since $\fs$ can be obtained as the sup
of functionals of the form \eqref{funcmeas} with $J$ finite,
we may assume that $J=\{x_1,\dots,x_m\}$.
Let $\mu_n \weakst \mu$ \wlystar in the sense of measures, and let $\lambda$
be the \wstar limit (up to a subsequence) of $|\mu_n \res J|$.
Let $J:=J_1 \cup J_2$, with
$$
J_1:=\{t \in J\,:\,|\mu(\{t\})| \ge \psi(t)\},
\quad\quad
J_2:=J \setminus J_1.
$$
Let $\eps>0$ be such that
$$
\bigcup_{x_i \in J_2} {\bar B}_\eps(x_i) \subseteq I
$$
and such that for all $n$
$$
|\mu_n| \left( \bigcup_{x_i \in J_2} \partial {\bar B}_\eps(x_i) \right)
=|\mu| \left( \bigcup_{x_i \in J_2} \partial {\bar B}_\eps(x_i) \right)=0.
$$
Let us set
$$
I_1:=I \setminus \bigcup_{x_i \in J_2} {\bar B}_\eps(x_i),
\quad\quad
I_2:=\bigcup_{x_i \in J_2} B_\eps(x_i).
$$
Let $\fs_1$ and $\fs_2$ denote the restriction of $\fs$ to $\ms_b(I_1;\R^k)$
and $\ms_b(I_2;\R^k)$
respectively. We have
$$
\liminf_n \fs(\mu_n) \ge \liminf_n \fs_1(\mu_n \res I_1)+
\liminf_n \fs_2(\mu_n \res I_2).
$$
We notice that
$$
\fs_1(\mu_n \res I_1) \ge \gs_1(\mu_n \res I_1)
$$
where
$$
\gs_1(\eta):=\int_{I_1} f_1(|\phi_{\eta}|)\,dx+
\sum_{t \in S_{\eta}} \varphi(|\eta(\{t\})|)+a|\eta^c|(I_1)
$$
for all $\eta \in \ms_b(I_1;\R^k)$. By \cite[Thorem 5.2]{AFP} we have that
$$
\gs_1(\mu \res I_1) \le \liminf_n \gs_1(\mu_n \res I_1),
$$
so that
$$
\fs_1(\mu \res I_1) =\gs_1(\mu \res I_1) \le \liminf_n \fs_1(\mu_n \res I_1).
$$
On the other hand, we have
$$
\fs_2(\mu_n \res I_2)=\gs_2(\mu_n \res {I_2 \setminus J_2})+
\sum_{t \in J_2}\varphi(|\mu_n(\{t\})| \lor \psi(t)),
$$
where
$$
\gs_2(\eta):=\int_{I_2} f_1(|\phi_{\eta}|)\,dx+
\sum_{t \in S_{\eta}} \varphi(|\eta(\{t\})|)+a|\eta^c|(I_2)
$$
for all $\eta \in \ms_b(I_2;\R^k)$. We have
\begin{align}
\liminf_n \fs_2(\mu_n \res I_2)
&\ge \gs_2(\mu \res {I_2 \setminus J_2})
+\sum_{t \in J_2}\varphi(\lambda(\{t\}) \lor \psi(t))
\\
\nonumber
&\ge \gs_2(\mu \res {I_2 \setminus J_2})
+\sum_{t \in J_2}\varphi(\psi(t)).
\end{align}
We deduce
$$
\fs_2(\mu \res I_2) =\gs_2(\mu \res {I_2 \setminus J_2})+
\sum_{t \in J_2}\varphi(\psi(t)) \le
\liminf_n \fs_2(\mu_n \res I_2),
$$
and so we get
$$
\fs(\mu)=\fs_1(\mu \res I_1)+\fs_2(\mu \res I_2)
\le \liminf_n \fs(\mu_n).
$$
The proof is now concluded.
\end{proof}

\begin{lemma}
\label{lsc}
Let $\bar F\,:\,BV(\Om) \to [0,+\infty]$ be defined by
$$
\overline{F}(u):=
\int_\Om f_1(|\nabla u|) \,dx+
\int_{\Sg{g}{u} \cup \Gamma} \varphi(|[u]| \lor \psi)\,d\hs^{N-1}
+a|D^cu|,
$$
with $a=\varphi'(0)$ and $f_1$ as in \eqref{bulkrelax}. Then $\bar F$
is lower semicontinuous with respect to the \wstar \linebreak
topology of $BV(\Om)$.
\end{lemma}

\begin{proof}
We may assume without loss of generality that
\begin{equation}
\label{intfinite}
\hs^{N-1}(\Gamma)<+\infty,
\quad\quad
\int_\Gamma \varphi(\psi)\,d\hs^{N-1}<+\infty.
\end{equation}
Firstly consider the case $\varphi(s)>\eps s$ for some $\eps>0$
and $s \in [0,+\infty[$.
Following \cite[Theorem 5.4]{AFP}, we use Lemma \ref{lscmeasure} to
obtain the lower semicontinuity in the one dimensional case, and we
recover the $N$-dimensional case using a slicing argument.
\par
Let us consider $\Omt$ open and bounded in $\R^N$ such that
$\Omb \subset \Omt$,
and let us set $\Omp:=\Omt \setminus \partial_N \Om$.
The lower semicontinuity of ${\bar F}$ is
equivalent to the lower semicontinuity of
\begin{equation}
\label{ftilde}
F'(u):=\int_{\Omp} f_1(|\nabla u|)\,dx+
\int_{S(u) \cup \Gamma} \varphi(|[u]|\lor \psi)\,d\hs^{N-1}
+a|D^cu|(\tilde{\Om})
\end{equation}
defined for all $u \in BV(\Omp)$ with $u=g$ on $\Omp \setminus \Om$.
In order to prove the lower semicontinuity of $F'$,
it is convenient to localize the functional,
i.e. for every open set $A \subseteq \Omp$, and
for every $u \in BV(\Omp)$
with $u=g$ on $\Omp \setminus \Om$, we consider
\begin{equation}
\label{ftildeloc}
F'_A(u):=\int_A f_1(|\nabla u|)\,dx+
\int_{A \cap (S(u) \cup \Gamma)} \varphi(|[u]| \lor \psi)\,d\hs^{N-1}
+a|D^cu|(A).
\end{equation}
Let $e \in \R^N$ with $|e|=1$: for every open set $A \subseteq \Omp$,
and for every $u \in BV(\Omp)$ such that $u=g$ on $\Omp \setminus \Om$
let us set
\begin{equation}
\label{ftildeslice}
F'_{A,e}(u):=\int_A f_1(|\langle \nabla u, e \rangle|)\,dx+
\int_{A \cap (S(u) \cup \Gamma)} |\langle \nu, e \rangle|\,
\varphi(|[u]| \lor \psi)\,d\hs^{N-1}
+a|\langle D^cu, e \rangle|(A).
\end{equation}
Here $\nu$ denotes the normal to the rectifiable set
$S(u) \cup \Gamma$.
\par
By the general theory of slicing (see \cite[Section 3.11]{AFP}),
we have that
\begin{equation}
\label{fdecom}
F'_{A,e}(u)=\int_{\pi_e}G_{A^e_y}(u^e_y)\,dy,
\end{equation}
where $\pi_e$ is the hyperplane through the origin orthogonal to $e$,
$A^e_y:=A \cap \{y+te\,,\,t \in \R\}$,
$u^e_y(t):=u(y+te)$, and
\begin{equation}
\label{gslice}
G_{A^e_y}(u^e_y):=
\int_{A^e_y} f_1(|(u^e_y)'|)\,dt+
\sum_{t \in (A^e_y \cap (S(u^e_y) \cup \Gamma^e_y))}
\varphi(|[u]|(y+te) \lor \psi(y+te))
+a|D^cu^e_y|(A^e_y).
\end{equation}
Let us consider $u_n \in BV(\Omp)$ such that $u_n=g$ on
$\Omp \setminus \Om$, and such
that $u_n \weakst u$ \wlystar in $BV(\Omp)$.
Then up to a subsequence for a.e. $y \in \pi_e$, we have
$(u_n)^e_y \to u^e_y$
strongly in $L^1((\Omp)^e_y)$. For every open set
$A \subseteq \Omp$, we claim that
for a.e. $y \in \pi_e$ we have
\begin{equation}
\label{lscgslice}
G_{A^e_y}(u^e_y) \le \liminf_n G_{A^e_y}((u_n)^e_y).
\end{equation}
In fact, if $\liminf_n G_{A^e_y}((u_n)^e_y)<+\infty$, we get
$(u_n)^e_y \weakst u^e_y$ \wlystar in $BV((\Omp)^e_y)$.
Moreover notice that for a.e. $y \in \pi_e$
$(\partial_D \Om)^e_y$ is finite,
and by \eqref{intfinite} $\Gamma^e_y$ is countable.
So in relation with these $y$'s, it is sufficient to set
$I:=A^e_y$ and $J:=(\partial_D \Om)^e_y \cup \Gamma^e_y$, and
to apply Lemma \ref{lscmeasure} to
$\mu_n:=D(u_n)^e_y$.
\par
Using Fatou's Lemma, by \eqref{fdecom} and \eqref{lscgslice}
we get
\begin{equation}
F'_{A,e}(u) \le \liminf_n F'_{A,e}(u_n).
\end{equation}
Let $\lambda:=\leb^N+\hs^{N-1} \res (S(u) \cup \Gamma)+|D^cu|$,
and let $B:=\{e_j\}$ be countable and dense in
$S^{N-1}:=\{x \in \R^N\,:\,|x|=1\}$.
If $E \subseteq \Omp \setminus (S(u) \cup \Gamma)$ is a
$\leb^N$ negligible set
on which $D^cu$ is concentrated, let us define
$$
f'_j(x):=
\begin{cases}
f_1(|\langle \nabla u(x), e_j \rangle|)
&\text{ if }x \in \Omp \setminus (E \cup S(u) \cup \Gamma) \\ \\
|\langle \nu(x), e_j \rangle|\, \varphi(|[u](x)| \lor \psi(x))
& \text{ if }x \in S(u) \cup \Gamma \\ \\
a\frac{|\langle D^cu, e_j \rangle|}{|D^cu|}(x)
&\text{ if }x \in E,
\end{cases}
$$
where $\nu(x)$ denotes the normal to the rectifiable set
$S(u) \cup \Gamma$ at the point $x$.
\par
For every $A_1,\dots,A_k$ disjoint open subsets of
$\tilde{\Om}$, and for every $e_1,\dots,e_k \in B$,
since $f_1$ is increasing, we obtain that
\begin{align}
\label{lscquasi}
\liminf_n F'(u_n) &\ge
\sum_{j=1}^k \liminf_n F'_{A_j}(u_n)
\ge \sum_{j=1}^k \liminf_n F'_{A_j,e_j}(u_n) \\
&\ge \sum_{j=1}^k F'_{A_j,e_j}(u)=
\sum_{j=1}^k \int_{A_j} f'_j \,d\lambda.
\end{align}
Applying \cite[Lemma 2.35]{AFP}, we deduce that
\begin{equation}
\liminf_n F'(u_n) \ge \int_{A} \sup_j f'_j \,d\lambda=F'(u),
\end{equation}
so that the Lemma is proved under the assumption
$\varphi(s) \ge \eps s$.
\par
The general case
follows observing that setting $\varphi_\eps(s):=\varphi(s)+\eps s$,
and letting
$\overline{F}_{\eps}$ be the functional defined in
\eqref{fbar} with
$\varphi_\eps$ in place of $\varphi$, we have
$$
\overline{F}(u) \le \overline{F}_{\eps}(u) \le
\overline{F}(u)+2\eps|Du|(\Om).
$$
\end{proof}

Let us now come to the proof of Proposition \ref{relaxation}.

\begin{proof}[Proof of Proposition \ref{relaxation}]
We can assume without loss of generality that
\begin{equation}
\label{intfinite2}
\int_\Gamma \varphi(\psi)\,d\hs^{N-1}<+\infty.
\end{equation}
Following Lemma \ref{lsc}, let us consider
$\Omt$ open and bounded in $\R^N$ such that $\Omb \subset \Omt$,
and let us set $\Omp:=\Omt \setminus \partial_N \Om$.
Let us consider the functional
\begin{equation}
\label{ftildetoberelax}
F'(u):=
\begin{cases}
\int_{\Om}f(|\nabla u|)\,dx+
\int_{S(u) \cup \Gamma} \varphi(|[u]|\lor \psi)\,d\hs^{N-1}
& \text{if }u \in SBV(\Omp), u=g \text{ on } \Omp \setminus \Om
\\ \\
+\infty &\text{otherwise in }BV(\Omp).
\end{cases}
\end{equation}
The relaxation result of Proposition \ref{relaxation} is
equivalent to prove that the relaxation
of \eqref{ftildetoberelax} under the \wstar topology of
$BV(\Omp)$ is
\begin{equation}
\label{ftilderelax}
\overline{F'}(u):=
\int_{\Om}f_1(|\nabla u|)\,dx+
\int_{S(u) \cup \Gamma} \varphi(|[u]|\lor \psi)\,d\hs^{N-1}+
a|D^cu|(\Omp)
\end{equation}
if $u \in BV(\Omp)$, $u=g$ on $\Omp \setminus \Om$, and
$\overline{F'}(u)=+\infty$
otherwise in $BV(\Omp)$.
\par
Following \cite{BBB}, it is useful to introduce the
localized version of \eqref{ftildetoberelax}; namely
for all open set $A \subseteq \Omp$ let us set
\begin{equation}
\label{flocal}
F'(u,A):=
\int_{A \cap \Om} f(|\nabla u|)\,dx+
\int_{A \cap (S(u) \cup \Gamma)}
\varphi(|[u]| \lor \psi)\,d\hs^{N-1}
\end{equation}
if $u \in SBV(\Omp)$, $u=g$ on $\Omp \setminus \Om$, and
$F'(u,A)=+\infty$ otherwise in $BV(\Omp)$. Let us indicate by
$\overline{F'}(u,A)$ the relaxation
of \eqref{flocal} under the \wstar topology of $BV(\Omp)$.
\par
Arguing as in \cite[Proposition 3.3]{BBB}, we have that for every
$u \in BV(\Omp)$, $\overline{F'}(u,\cdot)$ is
the restriction to the family $\as(\Omp)$ of all open subsets of
$\Omp$ of a regular Borel measure.
Since for all $u \in SBV(\Om')$ with $u=g$ on $\Omp \setminus \Om$
and for all $A \in \as(\Omp)$
we have
\begin{multline}
\int_{A \cap \Om} f(|\nabla u|)\,dx+
\int_{A \cap S(u)} \varphi(|[u]|)\,d\hs^{N-1}
\le
F'(u,A) \\
\le \int_{A \cap \Om} f(|\nabla u|)\,dx+
\int_{A \cap S(u)} \varphi(|[u]|)\,d\hs^{N-1}
+\int_{A \cap \Gamma} \varphi(\psi)\,d\hs^{N-1},
\end{multline}
by \cite[Theorem 3.1]{BBB} we obtain that for all
$u \in BV(\Om')$ with $u=g$ on $\Omp \setminus \Om$
and for all $A \in \as(\Omp)$ with $A \cap \partial_D \Om =\emptyset$
\begin{multline}
\label{comparison}
\int_{A \cap \Om} f_1(|\nabla u|)\,dx+
\int_{A \cap S(u)} \varphi(|[u]|)\,d\hs^{N-1}+a|D^cu|(A)
\le
F'(u,A) \\ \le
\int_{A \cap \Om} f_1(|\nabla u|)\,dx+
\int_{A \cap S(u)} \varphi(|[u]|)\,d\hs^{N-1}+a|D^cu|(A)
+\int_{A \cap \Gamma} \varphi(\psi)\,d\hs^{N-1}.
\end{multline}
As a consequence of \eqref{comparison}, we deduce that
$$
\overline{F'}(u,\cdot) \res
(\Omp \setminus (S(u) \cup \Gamma \cup \partial_D \Om))=
f_1(|\nabla u|) \,d\leb^N \res \Om+a|D^cu|.
$$
In order to evaluate
$\overline{F'}(u,\cdot) \res (S(u) \cup \Gamma \cup \partial_D \Om)$,
we notice that for all $A \in \as(\Omp)$
and for all $u \in SBV(\Omp)$ with $u=g$ on $\Omp \setminus \Om$
$$
\int_{A \cap \Om} f_1(|\nabla u|)\,dx+
\int_{A \cap (S(u) \cup \Gamma)}
\varphi(|[u]| \lor \psi)\,d\hs^{N-1}+a|D^cu|(A)
\le
F'(u,A),
$$
and since the left hand side is lower semicontinuous by Lemma
\ref{lsc}, we get that
for all $u \in BV(\Omp)$ with $u=g$ on $\Omp \setminus \Om$
\begin{equation}
\label{comparison2}
\int_{A \cap \Om} f_1(|\nabla u|)\,dx
+\int_{A \cap (S(u) \cup \Gamma)}
\varphi(|[u]| \lor \psi)\,d\hs^{N-1}+a|D^cu|(A)
\le
\overline{F'}(u,A).
\end{equation}
By outer regularity of $\overline{F'}(u,\cdot)$
we conclude that
$$
\overline{F'}(u,E) \ge
\int_{E} \varphi(|[u]| \lor \psi)\,d\hs^{N-1}
$$
for all Borel sets
$E \subseteq S(u) \cup \Gamma \cup \partial_D \Om$.
We have to prove the opposite inequality, namely
$$
\overline{F'}(u,E) \le
\int_{E} \varphi(|[u]| \lor \psi)\,d\hs^{N-1}
$$
for all Borel sets $E \subseteq S(u) \cup \Gamma \cup \partial_D \Om$.
Without loss of generality we may assume that
$$
\int_{S(u)}\varphi(|[u]|)\,d\hs^{N-1}<+\infty,
$$
and by a truncation argument, we can suppose that
$u_{|\Om} \in L^\infty(\Om)$. Let $K$ be a compact subset of
$S(u) \cup \Gamma \cup \partial_D \Om$, $\eps>0$, and
let $A_\eps$ be open with $K \subseteq A_\eps$ and
$$
|Du|(A_\eps \setminus K)<\eps,
\quad\quad
\int_{(A_\eps \setminus K) \cap \Gamma}
\varphi(\psi)\,d\hs^{N-1}<\eps.
$$
We can find $u_h \in BV(\Omp)$ with
$u_h=g$ on $\Omp \setminus \Om$ and such that
$u_h$ is piecewise constant in $\Om$ (that is
$(u_h)_{|\Om} \in SBV(\Om)$ with $\nabla u_h=0$
in $\Om$),
$u_h \to u$ strongly in $L^\infty(\Om)$, and
$|Du_h|(A_\eps \setminus K)<\eps$.
Since $u_h$ is piecewise constant in $\Om$
we have for all $h$
\begin{equation}
\overline{F'}(u_h,A_\eps) \le
\int_{A_\eps \cap (S(u_h) \cup \Gamma)}
\varphi(|[u_h]| \lor \psi)\,d\hs^{N-1}.
\end{equation}
We conclude
\begin{multline}
\overline{F'}(u,A_\eps) \le \liminf_h
\overline{F'}(u_h,A_\eps) \le
\liminf_h \int_{A_\eps \cap (S(u_h) \cup \Gamma)}
\varphi(|[u_h]|\lor \psi)\,d\hs^{N-1} \\ \le
\int_{K \cap (S(u) \cup \Gamma)}
\varphi(|[u]| \lor \psi)\,d\hs^{N-1}+
a|Du_h|(A_\eps \setminus K)
+\int_{(A_\eps \setminus K) \cap \Gamma}
\varphi(\psi)\,d\hs^{N-1} \\
\le
\int_{K \cap (S(u) \cup \Gamma)}
\varphi(|[u]| \lor \psi)\,d\hs^{N-1}+(a+1)\eps
\end{multline}
so that, letting $\eps \to 0$ we obtain
$$
\overline{F'}(u,K) \le
\int_{K \cap (S(u) \cup \Gamma)}
\varphi(|[u]| \lor \psi)\,d\hs^{N-1}.
$$
Since $K$ is arbitrary in
$S(u) \cup \Gamma \cup \partial_D \Om$,
the proof is concluded.
\end{proof}

\section{Two auxiliary propositions}
\label{auxsec}
In this section we prove two propositions we used
in the proofs of the main theorems of the paper.
For all $h \in \N$ let us consider $f_h:\R^N \to [0,+\infty[$
such that for all
$\xi \in \R^N$
\begin{equation}
\label{coercivity}
f_h(\xi) \nearrow |\xi|^2,
\quad\quad
f_h(\xi) \ge \min \{|\xi|^2-1, b_h|\xi|\}
\end{equation}
with $b_h \to +\infty$, and let $\varphi_h:[0,+\infty[ \to [0,1]$
be such that
for all $s \in [0,+\infty[$
\begin{equation}
\label{phiabovegzerobis}
\varphi_h(s) \ge \min\{c_hs,d_h\}
\end{equation}
with $c_h \to +\infty$ and $d_h \nearrow 1$ for $h \to +\infty$.
\par
The following result holds.

\begin{proposition}
\label{basicgconv}
Let $\Om \subseteq \R^N$ be open and bounded, and let us consider
the functionals
$$
F_h(u):=
\sum_{i=1}^m\int_{\Om} f_h(\nabla u_i) \,dx+
\int_{S(u)}\varphi_h(|[u_1]| \lor \ldots \lor |[u_m]|)
\,d\hs^{N-1}+a_h|D^cu|(\Om)
$$
if $u=(u_1,\dots,u_m) \in BV(\Om;\R^m)$, and $F_h(u)=+\infty$
otherwise in $BV(\Om;\R^m)$.
Let $a_h \to +\infty$ for $h \to +\infty$.
If $F_h(u^h) \le C$ and $u^h \weakst u$ \wlystar in $BV(\Om;\R^m)$,
we have $u \in SBV(\Om;\R^m)$,
\begin{equation}
\label{wgrad}
\nabla u^h \weak \nabla u
\quad\quad \text{weakly in }L^1(\Om;\R^{m \times N}),
\end{equation}
\begin{equation}
\label{basicgradlsc}
\|\nabla u_i\|^2 \le \liminf_h \int_{\Om}f_h(\nabla u^h_i) \,dx,
\quad\quad i=1,\dots,m,
\end{equation}
and
\begin{equation}
\label{basiclsc}
\hs^{N-1}(S(u)) \le \liminf_h \int_{S(u^h)}
\varphi_h(|[u^h_1]| \lor \ldots \lor |[u^h_m]|)\,d\hs^{N-1}.
\end{equation}
\end{proposition}

\begin{proof}
Let us consider $u^h \in BV(\Om;\R^m)$ such that $F_h(u^h) \le C$ and
$u^h \weakst u$ \wlystar in $BV(\Om;\R^m)$. Notice that $(\nabla u^h)$
is equintegrable.
In fact if $r_h$ is such that for all $|\xi| \le r_h$
$$
|\xi|^2-1 \le b_h|\xi|
$$
we get for all $i=1,\dots,m$ and for all $E \subseteq \Om$
\begin{align}
\label{estgradg}
\int_E |\nabla u^h_i|\,dx& \le  \int_{\{|\nabla u^h_i|\le r_h\} \cap E}
|\nabla u^h_i|\,dx
+\int_{\{|\nabla u^h_i| > r_h\} \cap E} |\nabla u^h_i|\,dx \\
\nonumber
&\le \left( \int_{\{|\nabla u^h_i|\le r_h\} \cap E}
|\nabla u^h_i|^2\,dx \right)^{\frac{1}{2}}|E|^{\frac{1}{2}}+
\int_{\{|\nabla u^h_i| > r_h\} \cap E} |\nabla u^h_i|\,dx \\
\nonumber
&\le \left( \int_{\Om}
(f_h(\nabla u^h_i)+1)\,dx \right)^{\frac{1}{2}}|E|^{\frac{1}{2}}+
\frac{1}{b_h}\int_{\Om} f_h(\nabla u^h_i)\,dx \le \sqrt{C+1}|E|+
\frac{C}{b_h},
\end{align}
where $|E|$ denotes the Lebesgue measure of $E$. This proves that
$\nabla u_h$ is equintegrable. Up to a subsequence we may suppose
that for all $i=1,\dots,m$ we have
$$
\nabla u^h_i \weak g_i
\quad\quad
\text{weakly in }L^1(\Om;\R^N).
$$
Since $a_h \to +\infty$, we get $D^cu^h \to 0$ strongly in the sense
of measures.
\par
Let us consider $\psi:\R \to \R$ bounded and Lipschitz, and for all
$i=1,\dots,m$ let us consider the measures
$$
\mu^h_i(B):=D\psi(u^h_i)(B)- \int_B \psi'(u^h_i) \nabla u^h_i \,dx,
\quad\quad
\lambda^h_i(B):= \int_{S(u^h_i) \cap B} \varphi_h(|[u^h_i]|)\,d\hs^{N-1},
$$
where $B$ is a Borel set in $\Om$.
We have
\begin{equation}
\label{chain}
|D\psi(u^h_i)- \psi'(u^h_i) \nabla u^h_i \,d\leb^N| \le
||\psi||_{\varphi_h}\lambda^h_i+ ||\psi'||_{\infty} |D^cu^h_i|
\end{equation}
where
$$
||\psi||_{\varphi_h}:= \sup \left\{
\frac{\psi(t)-\psi(s)}{\varphi_h(|t-s|)}\,:\,t \not= s \right\}.
$$
Up to a subsequence we have
$$
\mu^h_i \weakst D\psi(u_i)-\psi'(u)g_i \,d\leb^N,
\quad\quad
\lambda^h_i \weakst \lambda_i
$$
\wlystar in the sense of measures, and so by \eqref{chain} and since
$D^cu^h \to 0$ strongly
in the sense of measures we get
$$
|D\psi(u_i)-\psi'(u)g_i \,d\leb^N| \le (\sup \psi -\inf\psi) \lambda_i.
$$
As a consequence of $SBV$ characterization
(see \cite[Proposition 4.12]{AFP}), we get that
$u_i \in SBV(\Om)$, $\nabla u_i=g_i$ and
$\hs^{N-1} \res S(u_i) \le \lambda_i$ for all $i=1,\dots,m$.
We deduce that \eqref{wgrad} holds.
\par
In order to prove \eqref{basicgradlsc}, for every $M>0$ let $g^M_i$
be the weak limit in $L^1(\Om)$
(up to a subsequence) of $|\nabla u^h_i| \wedge M$.
Since $f_h(\xi) \to |\xi|^2$ uniformly on $[0,M]$, we have
\begin{equation}
\label{gm}
\|g^M_i\|^2 \le \liminf_h \int_\Om f_h(\nabla u^h_i)\,dx.
\end{equation}
Then as $M \to +\infty$ and summing over $i$ we have that
\eqref{basicgradlsc} holds.
\par
Let us come to \eqref{basiclsc}. If $\lambda$ is the weak limit in
the sense of measures
of
$$
\lambda^h(A):= \int_{S(u^h) \cap A} \varphi_h(|[u^h_1]|
\lor \ldots \lor |[u^h_m]|)\,d\hs^{N-1},
$$
we have that $\lambda_i \le \lambda$ for all $i=1,\dots,m$.
Since we have $\hs^{N-1} \res S(u_i) \le
\lambda_i$ for all $i=1,\dots,m$, we deduce that
$\hs^{N-1} \res S(u) \le \lambda$. Then we get that
\eqref{basiclsc} holds.
\end{proof}

Let $a_h \to +\infty$, and let us consider
$\varphi_h:[0,+\infty[ \to [0,1]$ increasing, concave,
$\varphi(0)=0$, $\varphi'(0)=a_h$ and such that
\begin{equation}
\label{varphihto1}
\varphi_h \nearrow 1.
\end{equation}
For all $h \in \N$ let us consider $f_h:\R^N \to [0,+\infty[$ such
that for all $\xi \in \R$
\begin{equation}
\label{fhgzero}
f_h(\xi)=
\begin{cases}
|\xi|^2 &|\xi| \le \frac{a_h}{2} \\
\frac{a_h^2}{4}+a_h(|\xi|-\frac{a_h}{2}) &|\xi| \ge \frac{a_h}{2}.
\end{cases}
\end{equation}
Let us consider $\Om$ open and bounded in $\R^N$,
and let $\partial_D \Om \subseteq \partial \Om$ be
open in the relative topology,
$\partial_N \Om:=\partial \Om \setminus \partial_D \Om$.
Let us consider
$\bar \Gamma \tsub \Om$
with $\hs^{N-1}(\bar \Gamma)<+\infty$, and let $\bar \gamma$
be a positive function on $\bar \Gamma$.
Let us extend $\bar \gamma$ to $\Omb$ setting
$\bar \gamma=0$ outside $\bar \Gamma$.
Then the following proposition holds.

\begin{proposition}
\label{gammaconvzero}
Let $g \in H^1(\Om)$ with $\|g\|_\infty <+\infty$, and let us
consider the functionals
\begin{equation}
\label{funczero}
\fs_h(u):=
\int_{\Om} f_h(\nabla u)\,dx+\int_{\bar \Gamma \cup \Sg{g}{u}}
\varphi_h(|[u]| \lor \bar \gamma)\,d\hs^{N-1}
+a_h|D^cu|(\Om)
\end{equation}
if $u \in BV(\Om)$, where $f_h$ and $\varphi_h$ are as
in \eqref{fhgzero} and \eqref{varphihto1}.
Then
$$
\inf_{u \in BV(\Om)} \fs_h(u) \to \inf_{u \in BV(\Om)} \fs(u)
$$
where
\begin{equation}
\label{funczerolim}
\fs(u):=
\begin{cases}
\|\nabla u\|^2+\hs^{N-1}(\bar \Gamma \cup \Sg{g}{u})
&u \in SBV(\Om) \\
+\infty &u \in BV(\Om) \setminus SBV(\Om).
\end{cases}
\end{equation}
\end{proposition}

\begin{proof}
In order to deal with $\Sg{g}{u}$ as an internal jump,
let us consider $\tilde{\Om} \subseteq \R^N$ open and bounded, and
such that $\Omb \subseteq \tilde{\Om}$.
Let us set $\Om':=\tilde{\Om} \setminus \partial_N \Om$.
Let us consider the functionals
$$
\fs'_h(u):=\int_{\Om'} f_h(\nabla u)\,dx+
\int_{\bar \Gamma \cup S(u)} \varphi_h(|[u]| \lor \bar \gamma)\,d\hs^{N-1}
+a_h|D^cu|(\Om')
$$
for $u \in BV(\Om')$ such that $u=g$ on $\Om' \setminus \Om$,
and $\fs'_h(u)=+\infty$ otherwise in
$BV(\Om')$.
In order to prove the lemma we can equivalently prove that
$$
\inf_{u \in BV(\Om')} \fs'_h(u) \to \inf_{u \in BV(\Om')} \fs'(u),
$$
where
\begin{equation}
\label{funczerolimprime}
\fs'(u):=
\begin{cases}
\|\nabla u\|^2+\hs^{N-1}(\bar \Gamma \cup S(u))
&u \in SBV(\Om'),\, u=g \text{ on }\Om' \setminus \Om \\
+\infty &\text{otherwise in }BV(\Om').
\end{cases}
\end{equation}
Let us consider
\begin{equation}
\label{defgshp}
\gs'_h(u):=
\begin{cases}
\fs'_h(u) &u \in BV(\Om'),\, \|u\|_\infty \le \|g\|_\infty \\
+\infty &\text{otherwise in }BV(\Om').
\end{cases}
\end{equation}
We have that $(\gs'_h)_{h \in \N}$ is an increasing
sequence of functionals
which converges pointwise to
\begin{equation}
\label{defgsp}
\gs'(u):=
\begin{cases}
\fs'(u) &u \in BV(\Om'),\, \|u\|_\infty \le \|g\|_\infty \\
+\infty &\text{otherwise in }BV(\Om').
\end{cases}
\end{equation}
By Lemma \ref{lsc} $\gs'_h$ is lower semicontinuous 
with respect
to the \wstar convergence in $BV(\Om')$, and
so by \cite[Proposition 5.4]{dm} we deduce that $\gs'_h$
$\Gamma${-}converges to $\gs'$ in the
\wstar topology of $BV(\Om')$. As a consequence 
of Proposition
\ref{Gamma-conv-prop}, we deduce that
$\inf_{BV(\Om')}\gs'_h \to \inf_{BV(\Om')}\gs'$.
By a truncation argument, we have that
$\inf_{BV(\Om')}\gs'_h =\inf_{BV(\Om')}\fs'_h$
and $\inf_{BV(\Om')}\gs'=\inf_{BV(\Om')}\fs'$, 
so that the proposition
is proved.
\end{proof}

\bigskip
\bigskip
\centerline{ACKNOWLEDGMENTS}
\bigskip\noindent
The author wishes to thank Gianni Dal Maso
for many helpful and interesting discussions.


\end{document}